\documentclass{amsart}
   \usepackage{graphicx}
\usepackage{amscd,amsmath,amsopn,amssymb,amsthm,multicol}
\usepackage{hyperref}
\usepackage[all]{xy}
\usepackage{amscd}
\usepackage{lscape}
\usepackage{slashed}
\usepackage{graphicx}
\usepackage{color}
\usepackage{lscape}
 \usepackage{cancel}
   \usepackage{harmony}
   \usepackage{comment}
   \usepackage{wasysym}
\usepackage{mathrsfs}
\DeclareMathAlphabet{\mathscrbf}{OMS}{mdugm}{b}{n}
\numberwithin{equation}{section}
\definecolor{dblue}{rgb}{0.01,0.01,0.44}
\definecolor{red}{rgb}{0.57,0.11,0.15}

 \input ulem.sty
 
\DeclareMathOperator{\Ad}{Ad}

\DeclareMathOperator{\Id}{Id}

\DeclareMathOperator{\Cl}{C\ell}
\DeclareMathOperator{\Ric}{Ric}
\DeclareMathOperator{\Sca}{Scal}

\DeclareMathOperator{\ke}{Ker}

\DeclareMathOperator{\Scho}{Sch}
\DeclareMathOperator{\Spin}{Spin}

\newcommand{\fr}{\mathfrak}
\newcommand{\al}{\alpha}
\newcommand{\be}{\beta}

\newcommand{\bb}{\mathbb}
\newcommand{\cal}{\mathcal}

\DeclareMathOperator{\SO}{SO}

 \DeclareMathOperator{\SU}{SU}
  \DeclareMathOperator{\Spec}{Spec}

\DeclareMathOperator{\G}{G}

\DeclareMathOperator{\Ss}{S}

   \newtheorem{lemma} {Lemma} [section]
\newtheorem{theorem}[lemma]{Theorem}
\newtheorem{remark}[lemma] {Remark}
\newtheorem{prop} [lemma]{Proposition}

\newtheorem{corol}[lemma] {Corollary}
\newtheorem{example}[lemma] {Example}

\begin{document}   

 \title
{A new $\frac{1}{2}$-Ricci type formula on the spinor bundle and applications} % which induces the generalized Schr\"odinger-Lichnerowicz  formula}
 \author{Ioannis Chrysikos}
  \address{Dipartimento di Matematica G. Peano, Universit\'a  di Torino, Via Carlo Alberto 10, 10123 Torino, Italy}
 \email{ioannis.chrysikos@unito.it}
% \curraddr{Department of Mathematics and Statistics,
%Case Western Reserve University, Cleveland, Ohio 43403}
 %\subjclass[2000]{Primary 54C40, 14E20; Secondary 46E25, 20C20}

%\date{January 1, 2001 and, in revised form, June 22, 2001.}


%\keywords{Differential geometry, algebraic geometry}

 % \medskip
%\noindent
%\thanks{The   author    was full-supported   
 % by Masaryk University under the Grant Agency of Czech Republic, project no.14-2464P}
 
%\thanks{The  author    was full-supported   
  %by Masaryk University under the Grant Agency of Czech Republic, project no. P 201/ 12/ G028}
 %\medskip

 \begin{abstract}
  Consider a Riemannian spin manifold $(M^{n}, g)$ $(n\geq 3)$ endowed with a non-trivial 3-form $T\in\Lambda^{3}T^{*}M$, such that $\nabla^{c}T=0$, where $\nabla^{c}:=\nabla^{g}+\frac{1}{2}T$  is the metric connection with skew-torsion $T$. In this note  we introduce a generalized  $\frac{1}{2}$-Ricci  type  formula for the spinorial action of the  Ricci endomorphism $\Ric^{s}(X)$, induced by the one-parameter family of metric connections $\nabla^{s}:=\nabla^{g}+2sT$.  This new  identity   extends   a result  described  by   Th.~Friedrich and E.~C.~Kim, about the action of the Riemannian Ricci endomorphism on spinor fields, and  allows us to present a series of applications.  For example, we describe  a new alternative proof of the generalized Schr\"odinger-Lichnerowicz formula related to the  square of the  Dirac operator $D^{s}$, induced by   $\nabla^{s}$, under the condition $\nabla^{c}T=0$.  In the same case, we   provide  integrability conditions for  $\nabla^{s}$-parallel spinors, $\nabla^{c}$-parallel spinors and     twistor spinors with torsion. We illustrate  our conclusions for some non-integrable structures satisfying our assumptions, e.g. Sasakian manifolds,  nearly K\"ahler manifolds  and nearly parallel $\G_2$-manifolds, in dimensions 5, 6 and 7, respectively.   %We show  that the restriction of the extended $\frac{1}{2}$-Ricci formula  to $\ker(\cal{P}^{s})$, coincides with the  twistorial $\frac{1}{2}$-Ricci type  formula  with respect to $\nabla^{s}$. This formula was recently introduced  by the author   in \cite[Lem.~2.2]{Chrysk2},  by a different   method. %We also examine $\nabla^{c}$-parallel spinors, and other special spinors.
   
 
   \end{abstract} 
\maketitle  
% \tableofcontents
 
\section{Introduction}\label{intro}

\noindent Let $(M^{n}, g)$ $(n\geq 3)$  be  a connected  Riemannian spin manifold  endowed with   a non-trivial 3-form $T\in\Lambda^{3}T^{*}M$. Consider the one-parameter family of connections $\{\nabla^{s} : s\in\bb{R}\}$, given by
\[
\nabla^{s}=\nabla^{g}+2sT.
\]
This is a line of  metric connections  with totally skew-symmetric torsion $T^{s}=4sT$, which  joins the connection $\nabla^{1/4}\equiv \nabla^{c}$ with torsion $T$, with the Levi-Civita connection $\nabla^{0}\equiv\nabla^{g}$.  By an  abuse of notation next we shall refer to $\nabla^{c}$  by the term ``characteristic connection''.   Let us denote by $\Ric^{s}$  the Ricci tensor  induced by $\nabla^{s}$.  In this note we    focus on  the action of the associated  Ricci endomorphsim $\Ric^{s}(X)$ $(X\in\Gamma(TM))$,  on the corresponding spinor bundle $\Sigma^{g}M$.  Under the condition $\nabla^{c}T=0$, and for any arbitrary spinor field $\varphi\in\cal{F}^{g}:=\Gamma(\Sigma^{g}M)$, we show that this action    can be described in terms of  the  Dirac operator  $D^{s}$ $(s\in\bb{R})$ induced by $\nabla^{s}$.  This takes place in   Section \ref{1/2}, where we provide  the following  {\it (generalized) $\frac{1}{2}$-Ricci type  formula}  (see Lemma \ref{LEMA})
\begin{eqnarray}
\frac{1}{2}\Ric^{s}(X)\cdot \varphi &=& D^{s}(\nabla^{s}_{X}\varphi)-\nabla^{s}_{X}(D^{s}\varphi)-\sum_{j=1}^{n}e_{j}\cdot \Big[\nabla^{s}_{\nabla^{s}_{e_{j}}X}\varphi+4s\nabla^{s}_{T(X, e_{j})}\varphi\Big]\nonumber\\
&&+s(3-4s)(X\lrcorner \sigma_{T})\cdot\varphi, \label{gen1/2ric}
\end{eqnarray}
 for any arbitrary  vector field $X\in\Gamma(TM)$,  spinor field $\varphi\in\cal{F}^{g}$  and $s\in\bb{R}$,   where $\sigma_{T}$  is the 4-form 
 \[
    \sigma_{T}:=\frac{1}{2}\sum_{i=1}^{n}(e_{i}\lrcorner T)\wedge (e_{i}\lrcorner T).
    \] 
\noindent   From now on we shall mainly refer to $( \ref{gen1/2ric})$  by  the term {\it $\frac{1}{2}$-$\Ric^{s}$-formula},  or {\it $\frac{1}{2}$-$\Ric^{s}$-identity}.  This  can be viewed as the analogue  of   the   {\it Riemannian   $\frac{1}{2}$-Ricci formula}, or in short {\it $\frac{1}{2}$-$\Ric^{g}$-formula}, introduced by  Friedrich and  Kim  in \cite[Lem.~1.2]{FKim}. The latter relates  the Ricci endomorphism of the Levi-Civita  connection with the Riemannian Dirac operator, i.e. \begin{equation}
\frac{1}{2}\Ric^{g}(X)\cdot \varphi =D^{g}(\nabla^{g}_{X}\varphi)-\nabla^{g}_{X}(D^{g}\varphi)-\sum_{j=1}^{n}e_{j}\cdot \nabla^{g}_{\nabla^{g}_{e_{j}}X}\varphi. \label{1/2FK} 
\end{equation}

   In \cite{FKim} it was shown   that {\it  the  $\frac{1}{2}$-$\Ric^{g}$-identity    is stronger than the Schr\"odinger-Lichnerowicz  formula  associated to  the Riemann Dirac operator} $D^{g}\equiv D^{0}$, in the sense that  the first formula induces the second one, after a contraction. %  (\ref{ricg}) induces (\ref{SL}).
Here, we extend this result  by proving  that  the new $\frac{1}{2}$-$\Ric^{s}$-formula  induces the  corresponding {\it generalized formula of  Schr\"odinger-Lichnerowicz type},  associated to the Dirac operator $D^{s}$   (see for example  \cite[Thm.~3.1]{FrIv},  \cite[Thm.~6.1]{AF}  or  \cite[Thm.~3.2]{Agr03}), under the condition $\nabla^{c}T=0$.  % the {\it generalized SL-formula} associated to $D^{s}$.
 Therefore, when  the torsion form $T$ is $\nabla^{c}$-parallel we provide a new   proof for  this fundamental formula which is  different than the traditional proofs, compare for instance with \cite{FrIv, Agr03, AF}. %  There,   the   generalized SL-formula occurs by computing  the square  $(D^{s})^{2}$;   in contrast, our proof relies on the new  $\frac{1}{2}$-$\Ric^{s}$-formula and a contraction. 
% Hence, for parallel skew-torsion  the   $\frac{1}{2}$-$\Ric^{s}$-formula (which for $s=0$ reduces to (\ref{1/2FK})), provides  us with all the necessary information in order to   verify  that the previous conclusion of \cite{FKim}, can be successively extended  to the present 1-parameter family   of metric connections.

%\vskip 0.2cm
%\noindent  % To be more precise, a   formula of  Schr\"odinger-Lichnerowicz type (in short, SL-type) for the square of $D^{c}$ and more general of $D^{s}$ is known by \cite[Thm.~3.1]{FrIv}, \cite[Thm.~3.2]{Agr03} and \cite[Thm.~6.1]{AF}. Similar formulas of SL-type for the so-called {\it cubic Dirac operator} $D^{1/12}:=D^{g}+\frac{1}{4}T$, were known before    by \cite[Thm.~2.16]{Kos}, or \cite[Thm.~3.3]{Agr03}, for example (see also  \cite[Thm.~2.3]{Bismut} and \cite[pp.~5-6]{Dalakov} for remarks on Bismut's paper and moreover the SL-formula of the cubic Dirac operator in dimension 4).  Important for our approach  is the case with parallel torsion,  $\nabla^{c}T=0$, which we recall in Section \ref{1/2}, see Thereom \ref{parallel}. 

  \vskip 0.1cm
  The new $\frac{1}{2}$-Ricci type identity, being stronger than the generalized SL-formula for $D^{s}$, has several nice applications. In fact,  it is a  spinorial identity which      reproduces  all  $\frac{1}{2}$-Ricci type formulas  associated to $\nabla^{s}$ (in the sense of \cite{FKim}), even for $s=0$, but also other known results.  For example, in \cite{Chrysk2} we have recently introduced  a   {\it twistorial $\frac{1}{2}$-$\Ric^{s}$-formula} for  {\it twistor spinors with torsion} with respect to the family $\nabla^{s}$. Such spinors are elements in the kernel of the Penrose operator $\cal{P}^{s}$, induced by $\nabla^{s}$. When $T$ is $\nabla^{c}$-parallel and $M^{n}$ is compact,  in \cite[Corol.~3.2]{ABK} it was shown that twistor spinors with torsion realize the equality case of  an  estimate for the first eigenvalue of the square of the cubic Dirac operator $D^{1/12}=D^{g}+\frac{1}{4}T$, under some additional geometric assumptions (e.g. constant scalar curvature).  
The twistorial $\frac{1}{2}$-$\Ric^{s}$-formula (\cite[Lem.~2.2]{Chrysk2})  appears also  under the condition $\nabla^{c}T=0$  and in the context of  spin geometry with (parallel) skew-torsion, it establishes   the analogue   of a basic result of Lichnerowicz \cite{Lic} (see also     \cite[p.~123]{Fr}   or \cite[Prop.~A.2.1.(3a)]{Ginoux}).  In   Section \ref{1/2} we  obtain  the  twistorial $\frac{1}{2}$-$\Ric^{s}$-formula via a new and easier method, in particular we prove that it coincides with the restriction of the $\frac{1}{2}$-$\Ric^{s}$-identity to   the kernel   of the   twistor operator $\cal{P}^{s}$ (see Theorem \ref{prof}). % c
 
 \vskip 0.1cm
 Next    we   proceed with an examination of    $\nabla^{c}$-parallel spinors and more general $\nabla^{s}$-parallel spinors.  Recall that  when  $\nabla^{c}$ is the characteristic connection of a non-integrable $G$-structure on $(M^{n}, g)$ (in terms for example of \cite{FrIv}), then the   condition $\nabla^{c}\varphi=0$ for some non-trivial spinor field $\varphi$, imposes restrictions to the holonomy group ${\rm Hol}(\nabla^{c})\subset G$. % since the spinor holonomy representation has to have a fixed point.   
Here, when $T$ is $\nabla^{c}$-parallel,  we deduce that a  non-trivial spinor field  $\varphi_{0}\in\cal{F}^{g}$  which is parallel with respect to  $\nabla^{s}$  for some parameter $s\in\bb{R}$, must satisfy the  following equations (for any $X\in\Gamma(TM)$ and for the same $s$)
 \[
   \Ric^{s}(X)\cdot\varphi_{0}=2s(3-4s)(X\lrcorner \sigma_{T})\cdot\varphi_{0}, \quad   \quad \Sca^{s}\cdot \ \varphi_{0}=-8s(3-4s)\sigma_{T}\cdot\varphi_{0}.
\]
  For $s=0$  this yields the well-known $\Ric^{g}$-flatness of $(M^{n}, g)$, while for $\nabla^{c}$-parallel spinors we obtain   the  conditions given by  Friedrich and  Ivanon \cite{FrIv}, i.e. $\ \Ric^{c}(X)\cdot\varphi_{0}=(X\lrcorner \sigma_{T})\cdot\varphi_{0}$ and $\sigma_{T}\cdot\varphi_{0}=-\frac{\Sca^{c}}{4}\cdot\varphi_{0}$.  Our  most  interesting  result  is  related with $\nabla^{c}$-parallel spinors. Such spinors   have applications in theoretical physics, especially in type II string theory, where basic models are described in terms of   a metric connection with skew-torsion and the corresponding parallel spinors represent  the preserved supersymmetries (for more background we refer to \cite{IvPap, FrIv, FrIv2}).  In  Section \ref{psps}  we  present    the explicit action of the endomorphism  $\Ric^{g}(X)$ and more general $\Ric^{s}(X)$  on ${\rm Ker}(\nabla^{c})$, where $\nabla^{c}$ is any metric connection with skew-torsion $T$ such that $\nabla^{c}T=0$ (this means, without assuming   that $\nabla^{c}$ is the characteristic connection of some underlying special structure).  In particular, we provide the following remarkable formula  (see   Theorem \ref{newpar} and Corollaries \ref{newpar2}, \ref{Ssym})
  \begin{eqnarray*}
  \Ric^{s}(X)\cdot\varphi_{0}&=&\frac{(16s^{2}-1)}{4} \sum_{j=1}^{n}T(X, e_{j})\cdot (e_{j}\lrcorner T)\cdot\varphi_{0}+\frac{(16s^2 + 3)}{4}(X\lrcorner \sigma_{T})\cdot\varphi_{0}\\
  &=&\Ric^{c}(X)\cdot\varphi_{0}-\frac{(16s^{2}-1)}{4} S(X)\cdot\varphi_{0},
  \end{eqnarray*}
 for any $\nabla^{c}$-parallel spinor $\varphi_{0}$ and $X\in\Gamma(TM)$, where the endomorphism $S(X)$ is given by
 \[
 S(X):=-X\lrcorner \sigma_{T} +\sum_{j=1}^{n}e_{j}\cdot\big(T(X, e_{j})\lrcorner T). %]-(X\lrcorner \sigma_{T}).
 \]
   Then, we specialise on some types of non-integrable geometric structures satisfying our assumptions, e.g. 5-dimensional Sasakian manifolds, 6-dimensional nearly K\"ahler manifolds and 7-dimensional nearly parallel $\G_2$-manifolds.   We    illustrate our  integrability conditions  and describe the action of $\Ric^{s}(X)$ on the corresponding $\nabla^{c}$-parallel spinors (adapted to the particular special structure). For the Sasakian case,  our result     (see Theorem \ref{genparsak})  nicely  extends the integrability conditions given in \cite[Thm.~7.3, 7.6]{FrIv2}, for any $s\in\bb{R}$. %Hence,   the new $\frac{1}{2}$-$\Ric^{s}$-formula essentially allows us to extend result by these cited works. 
   For   nearly K\"ahler structures  and   nearly parallel $\G_2$-structures we recover some of our conclusions   in \cite{Chrysk2},  however by a new method (see Corollary \ref{mine22}). 

 \vskip 0.1cm
  A final contribution of this  note is related with the following first-order differential operator acting on   spinors,  $\slashed{D}^{s}(\varphi):=\sum_{j}(e_{j}\lrcorner T)\cdot \nabla^{s}_{e_{j}}\varphi$. This operator  is  included in the expression of $(D^{s})^{2}$  and in the compact case can be viewed as the main obstruction to  a universal  estimate of the lowest eigenvalue of    $(D^{s})^{2}$, see   \cite{Agr03, FrIv, AF}. % With the aim to overpass this problem, one makes the rescaling trick $s\to s'=s/3$, which finally for $s=1/4$ induces the cubic Dirac operator $D^{1/12}$, see \cite{Agr03, AF}. %The cubic Dirac operator is the generalized Dirac operator associated to the connection  with torsion $T/3$, where $T$ is the torsion of $\nabla^{c}$, and indeed an  advantage of the associated SL-formula is that allows us to write down explicitly universal estimates  of the first eigenvalue of  $(D^{1/12})^{2}$, see for instance \cite[Cor.~3.1]{Agr03} or \cite[Thm.~2.3]{ABK}.  
In  Section \ref{hope2} we examine  $\slashed{D}^{s}$ and describe some special  kinds of $\slashed{D}^{s}$-eigenspinors (see Proposition \ref{parslash},  \ref{endsla}). We also provide   examples.  Different classes of $\slashed{D}^{s}$-eigenspinors will be presented in a forthcoming work. %see \cite{Agr03, AF}. We also examine under which conditions, eigenspinors of $\slashed{D}^{s}$ can induce special spinors with respect to $\nabla^{s}$ (see Theorem \ref{}, \ref{}) and we present some examples.

 \vskip 0.2cm
 \noindent {\bf Thanks.} 
 The  author is  a Marie Curie fellow of the Istituto Nazionale di Alta Matematica ({INdAM}).   He  thanks Prof. Anna Fino (Torino) for discussions and comments.  % and  acknowledges the Dipartimento di Matematica ``Giuseppe Peano'', Universit\'a  di Torino, for its hospitality. 
 
\section{Preliminaries}\label{spinin}
 \noindent With the aim  to set  up our conventions relevant to subsequent computations, we begin by  recalling  basic facts from spin geometry.  Since all this material is well-known, we  provide  an exposition  only of the most useful notions (without proofs)  and for any further and detailed information, we refer   to \cite{Fr, Ginoux, FrIv, AF, ABK, Agr03}.

\subsection{Spin geometry with torsion}\label{nablas}
Consider   an oriented connected  Riemannian  manifold $(M^{n}, g)$ $(n\geq 3)$ endowed with a spin structure, i.e. a  $\Spin(n)$-principal bundle $\widetilde{P}^{g}:=\widetilde{\SO}(M, g)\to M$  together with a 2-fold covering $\Lambda^{g} : \widetilde{P}^{g}  \to P^{g}$,  such that $\Lambda^{g}(\tilde{p}g) = \Lambda^{g}(\tilde{p})\Ad(g)$  for any $\tilde{p}\in \widetilde{P}^{g}$ and $g\in\Spin_{n}$. Here,   and for the following  of this article  we denote by $P^{g}:=\SO(M, g)$  the $\SO_{n}$-principal bundle of positively oriented orthonormal frames of $M$.  We also remind that for $n\geq 3$, the spin group $\Spin_{n}$ is  the universal covering of $\SO_{n}$ and  $\Ad : \Spin_{n}\to\SO_{n}$  denotes   the double covering map. 
Via the spin representation (which we agree to denote by $\kappa_{n}$), we  associate to $\widetilde{P}^{g}$    a complex vector bundle $\Sigma^{g}M\to M$,  the so-called   spinor bundle  $\Sigma^{g}M:=\widetilde{P}^{g}\times_{\kappa_{n}}\Delta_{n}=\widetilde{P}^{g}\times_{\Spin_{n}}\Delta_{n}$, where $\Delta_{n}$ is the spin module. Notice that the spinor bundle  cannot be defined independently of a (semi)-Riemannian metric, in particular the definition of spinor fields, i.e. sections of the spinor bundle, depends in general on $g$,  in contrast to tensors.   For the following we set $\cal{F}^{g}:=\Gamma(\Sigma^{g}M)$ and   recall that $\Sigma^{g}M$  is endowed with a  $\Spin_{n}$-invariant Hermitian  product $\langle \ , \  \rangle$,  defined fiberwise as the natural Hermitian scalar product  that admits $\Delta_{n}$.  Its real part $( \ , \ ):={\rm Re}\langle  \ , \ \rangle$  induces a positive definite inner product on $\Sigma^{g}M$.   The  Clifford multiplication is the  bundle morphism   $\mu : TM\otimes_{\bb{R}}\Sigma^{g}M\to\Sigma^{g}M$, defined by  $\mu(X\otimes \varphi):=\kappa_{n}(X)(\varphi)=X\cdot\varphi$ and it naturally extends to differential forms $\mu : \Lambda(M)\otimes_{\bb{R}}\Sigma^{g}M\to\Sigma^{g}M$, 
\[
\omega\cdot\varphi:=\sum_{1\leq i_{1}<\cdots < i_{p}\leq n}\omega_{i_{1}\cdots i_{p}}e_{i_{1}}\cdot \ldots \cdot e_{i_{p}}\cdot \varphi,
\]
 for any  $\omega=\sum_{1\leq i_{1}<\cdots < i_{p}\leq n}\omega_{i_{1}\cdots i_{p}}e_{i_{1}}\wedge \ldots \wedge e_{i_{p}}\in\Lambda^{p}T^{*}M$.    Given a vector field $X$, let us denote by $X^{\flat}$   the dual 1-form, i.e. $X^{\flat}(u)=g(X, u)$. Then, any $X, Y\in T_{x}M$,  $\omega\in\Lambda^{p}T^{*}_{x}M$ and    $\varphi, \psi\in\Delta_{n}$ satisfy the following very useful properties (and similar for sections)%inherits the relations of the Clifford algebra,
%\[
%X\cdot Y\cdot \varphi+Y\cdot X\varphi=-2g(X, Y)\varphi, 
%\]
%  and   satisfies the relations 
  \begin{equation}\label{clif1}
 \left\{
 \begin{tabular} {r l l || r l l}
$-2g_{x}(X, Y)\varphi$ & $=$ & $X\cdot Y\cdot \varphi+Y\cdot X\cdot \varphi$ & $ \langle X\cdot\varphi, \psi\rangle$ & $=$ & $-\langle \varphi, X\cdot\psi\rangle$ \\ 
$X\cdot\omega$ & $=$ & $X^{\flat}\wedge\omega-X\lrcorner \omega$ & $\langle \omega\cdot \varphi, \psi\rangle$ & $=$ & $(-1)^{p(p+1)/2}\langle\varphi, \omega\cdot\psi\rangle$ \\
$\omega\cdot X$ & $=$ & $(-1)^{p}(X^{\flat}\wedge\omega+X\lrcorner \omega)$ & $-2(X\lrcorner \omega)$ & $=$ & $X\cdot\omega-(-1)^{p}\omega\cdot X$.
 \end{tabular}\right.
 \end{equation}

 \vskip 0.2cm
 \noindent %Now, any  metric connection  on $TM$ (corresponding to some connection 1-form $\omega : TP^{g}\to \fr{so}(n)$), lifts to a metric     connection in the spinor bundle $\Sigma^{g}M$ (corresponding to a connection 1-form $\widetilde{\omega} : T\widetilde{P}^{g} \to\fr{spin}(n)$ such that $\ad\circ \, \widetilde{\omega}=\omega\circ  d\Lambda^{g}$, with $\ad=d\Ad$). In particular, from now on we 
 From now on, let us assume     that {\it $(M^{n}, g)$  is equipped with a non-trivial 3-form $T\in\Lambda^{3}T^{*}M$}. We consider the one-parameter family of  metric connections   $\{\nabla^{s} :  s\in\bb{R}\}$ with  skew-torsion $4sT$;  this is defined by  
 \[
 g(\nabla^{s}_{X}Y, Z)=g(\nabla^{g}_{X}Y, Z)+2sT(X, Y, Z),
\] 
for any $X, Y, Z\in\Gamma(TM)$. The family $\nabla^{s}$ lifts to a family of  metric connections  on  $\Sigma^{g}M$, say $\nabla^{s} : \Gamma(\Sigma^{g}M)\to \Gamma(T^{*}M\otimes\Sigma^{g}M)$ (we keep the same notation), which explicitly reads by $\nabla^{s}_{X}\varphi=\nabla^{g}_{X}\varphi+s(X\lrcorner T)\cdot\varphi$.  % the Levi-Civita connection  lifts to a spinorial Levi-Civita connection in  $\widetilde{P}^{g}$.  This  induces  a metric covariant derivative on $\Sigma^{g}M$, which we denote by the same notation, $\nabla^{g} : \cal{F}^{g}\to \Gamma(T^{*}M\otimes\Sigma^{g}M)$.  
In  terms of some local orthonormal frame $\{e_{i}\}$, it is $\nabla^{s}_{X}\varphi=\nabla^{g}_{X}\varphi+s\sum_{i<j}^{n}T(X, e_{i}, e_{j})e_{i}\cdot e_{j}\cdot\varphi$  and the metric compatibility has the form $X\langle \varphi, \psi\rangle=\langle \nabla^{s}_{X}\varphi, \psi\rangle+\langle \varphi, \nabla^{s}_{X}\psi\rangle$. We also remind the Liebniz rule, 
 \[
 \nabla^{s}_{X}(Y\cdot\varphi)=(\nabla^{s}_{X}Y)\cdot \varphi+Y\cdot \nabla_{X}^{s}\varphi,\quad \nabla^{s}_{X}(\omega\cdot\varphi)=(\nabla^{s}_{X}\omega)\cdot \varphi+\omega\cdot \nabla_{X}^{s}\varphi,
 \]
  for any $X, Y\in\Gamma(TM)$, $\omega\in\Lambda^{p}T^{*}M$ and $\varphi, \psi\in\cal{F}^{g}$. %In terms of  a local representation $[\tilde{p}, \tilde{\varphi}]$ of $\varphi\in\cal{F}^{g}$  where  $\tilde{p}\in\Gamma_{U}(\widetilde{P}^{g})$,   $U\subset M$ is an open subset and   $\tilde{\varphi} : U\subset M\to\Delta_{n}$ a function, it holds that
% \[
%  \nabla^{s}_{X}\varphi=\Big[\tilde{p}, d\tilde{\varphi}(X)+\frac{1}{2}\sum_{i<j}g(\nabla^{s}_{X}e_{i}, e_{j})e_{i}\cdot e_{j}\cdot\tilde{\varphi}\Big].
% \]

\vskip 0.2cm
\noindent The (spinorial) curvature operator $\cal{R}^{s}_{X, Y}:=[\nabla^{s}_{X}, \nabla^{s}_{Y}]-\nabla^{s}_{[X, Y]} : \cal{F}^{g}\to\cal{F}^{g}$ associated to the covariant derivative $\nabla^{s}$ on the spinor bundle, satisfies the relation  $
\cal{R}^{s}_{X, Y}\varphi=\frac{1}{2}R^{s}(X\wedge Y)\cdot\varphi$, 
 where $R^{s}$ is the  curvature operator on 2-forms, induced by $\nabla^{s}$ at the tangent bundle level. Locally, one has the relations
\[
R^{s}(e_{i}\wedge e_{j}):=\sum_{k<l}R^{s}_{ijkl}e_{k}\wedge e_{l}, \quad   \cal{R}^{s}_{X, Y}\varphi=\frac{1}{2}\sum_{i<j}g(R^{s}(X, Y)e_{i}, e_{j})e_{i}\cdot e_{j}\cdot\varphi. 
\]
%When $\nabla^{c}T=0$, the curvature tensor $R^{s}$ is symmetric $R^{s}(X, Y, Z, W)=R^{s}(Z, W, X, Y)$ and for $\cal{R}^{s}$ locally we can write
% \[
% \cal{R}^{s}_{X, Y}\varphi=\frac{1}{2}\sum_{i<j}g(R^{s}(X, Y)e_{i}, e_{j})e_{i}\cdot e_{j}\cdot\varphi. 
%\]
For $s=0$, one obtains the (spinorial) Riemannian curvature operator $\cal{R}^{g}_{X, Y}$ and then, the  first Bianchi identity associated to $\nabla^{g}$ yields the well-known relation
\begin{equation}\label{1/2ricg}
\frac{1}{2}\Ric^{g}(X)\cdot\varphi=-\sum_{i=1}^{n}e_{i}\cdot \cal{R}^{g}_{X, e_{i}}\varphi,
\end{equation}
where $\Ric^{g}(X)$ is the Riemannian  Ricci endomorphism,  i.e. $\Ric^{g}(X, Y)=g(\Ric^{g}(X), Y)$, for any $X, Y\in\Gamma(TM)$.   For our family $\nabla^{s}$, and under the assumption $\nabla^{c}T=0$,   the associated Ricci tensor $\Ric^{s}$  remains symmetric (see also Remark \ref{ricci} below) and the Ricci endomorphism $\Ric^{s}(X)$ $(X\in\Gamma(TM))$  acts on spinors by  the rule 
\[
\Ric^{s}(X)\cdot\varphi:=\sum_{i}g(\Ric^{s}(X), e_{i})e_{i}\cdot\varphi=\sum_{i}\Ric^{s}(X, e_{i}) e_{i}\cdot\varphi.
\] In the same case,  Becker-Bender  proved  that  the first Bianchi  identity associated to $\nabla^{s}$ (cf. \cite[Thm.~B.1]{ABK})  induces  an analogue of (\ref{1/2ricg}), namely: 
\begin{lemma}\textnormal{(\cite[Lem.~1.13]{Julia})} \label{JUL1}
 Under the assumption $\nabla^{c}T=0$,  the  following relation holds  
   \[
\frac{1}{2}\Ric^{s}(X)\cdot\varphi =- \sum_{i}e_{i}\cdot \cal{R}^{s}_{X, e_{i}}\varphi +s(3-4s)(X\lrcorner \sigma_{T})\cdot\varphi, 
 \]
  for any arbitrary   vector field $X\in\Gamma(TM)$, spinor field $\varphi\in\cal{F}^{g}$ and $s\in\bb{R}$. \end{lemma}
%Here, by  $\sigma_{T}$  we denote  the following 4-form on $M$,  
 %\[
% \sigma_{T}:=\frac{1}{2}\sum_{i=1}^{n}(e_{i}\lrcorner T)\wedge (e_{i}\lrcorner T).
 %\]
 %An equivalent definition is given by $\sigma_{T}(X, Y, Z, W)=\fr{S}^{X, Y, Z}g(T(X, Y), T(Z, W))$.
% \begin{remark}\label{st}
%\textnormal{When the condition $\nabla^{c}T=0$ holds,  one can prove that $\nabla^{c}\sigma_{T}=0$, $\delta^{s}T=0$ for any $s\in\bb{R}$ and $dT=2\sigma_{T}$, see \cite{AF, FrIv, ABK}.}    
   %\end{remark}

  \begin{remark}\label{ricci}\textnormal{Let us denote by $d^{s}$ and $\delta^{s}$ the differential and co-differential induced by $\nabla^{s}$, i.e.
   \[
   d^{s}\omega:=\sum_{j}e_{j}\wedge \nabla^{s}_{e_{j}}\omega, \quad \delta^{s}\omega:=- \sum_{j}e_{j}\lrcorner \nabla^{s}_{e_{j}}\omega,
   \]
   respectively.  For $s=0$, we set $d:=d^{g}\equiv d^{0}$ and $\delta:=\delta^{g}\equiv\delta^{0}$.   By \cite{FrIv, AF, ABK} it is known that for some 3-form $T\in\Lambda^{3}T^{*}M$, it is $\delta T=\delta^{s}T$ for any $s\in\bb{R}$.  Notice also that the 4-from $\sigma_{T}$ has a distinctive role in the theory of  metric connections with skew-torsion. For example,   by using \cite[Lem.~2.4]{Agr03}  one deduces that $\sigma_{T}$ measures the difference of the differentials $d^{s}$ and $d$, i.e. $d^{s}T=dT-8s\sigma_{T}$. When the additional condition $\nabla^{c}T=0$ holds,  it is known  that $dT=2\sigma_{T}$ \cite{IvPap, FrIv} and $\nabla^{c}\sigma_{T}=0=\delta^{s}T$, for any $s\in\bb{R}$  \cite{AF, ABK}. In the same case, it is also true that  the curvature tensor $R^{s}$ is symmetric, $R^{s}(X, Y, Z, W)=R^{s}(Z, W, X, Y)$  and the same holds for the associated Ricci tensor    $\Ric^{s}(X, Y):=\sum_{i}R^{s}(X, e_{i}, e_{i}, Y)$, i.e. $   \Ric^{s}(X, Y)=\Ric^{s}(Y, X),$    for any $X,Y\in\Gamma(TM)$ and $s\in\bb{R}$, see \cite[Thm.~B.1]{ABK}. }% Notice that when $\nabla^{c}T\neq 0$, the sum $\sum_{i<j}\Big[\Ric^{s}(e_{i}, e_{j})-\Ric^{s}(e_{j}, e_{i})\Big]$ is not trivial. }%   In particular,  one can prove that $ \Ric^{s}(X, Y)=\Ric^{g}(X, Y)-4s^{2}S(X, Y)$,  where $S$ is the symmetric tensor   defined by $S(X, Y):=\sum_{i}g(T(X, e_{i}), T(Y, e_{i}))$. }%For completeness we also recall that the scalar curvature satisfies  the relation $\Sca^{s}=\Sca^{g}-24s^{2}\|T\|^{2}$, where $\|T\|^{2}=\frac{1}{3}\sum_{i<j}g\big(T(e_{i}, e_{j}), T(e_{i}, e_{j})\big)$ denotes the normalized  length of the 3-form $T$ and $\Sca^{g}$ the Riemannian scalar curvature.%The  $\Spin_{n}$-representation $\bb{R}^{n}\otimes\Delta_{n}$ splits into two  irreducible parts, i.e.  $\bb{R}^{n}\otimes\Delta_{n}=\ker(\mu)\oplus\Delta_{n}$  and this  decomposition  extends at a  bundle level, $TM\otimes\Sigma^{g}M=\ker(\mu)\oplus \Sigma^{g}M$ and consequently between sections $\Gamma(TM\otimes \Sigma^{g}M)=\Gamma(\ker(\mu))\oplus\cal{F}^{g}$. Let us denote by $p : \bb{R}^{n}\otimes\Delta_{n}\to\ker(\mu)$ the universal projection. In terms of a local orthonormal frame   $\{e_{i}\}$,  it is $ X\otimes\varphi\mapsto X\otimes\varphi +\frac{1}{n}\sum_{i}e_{i}\otimes e_{i}\cdot X\cdot\varphi$. 

%}
  \end{remark}  
  
 % \subsection{Differential operators on spinors}
% \vskip 0.1cm
 \noindent  We pass now on differential operators acting on spinors.  The   {\it (generalized)  Dirac operator} is the first-order differential operator on  $\cal{F}^{g}$,  defined by 
	\[
	D^{s}:=\mu\circ\nabla^{s} : \cal{F}^{g}\overset{\nabla^{s}}{\to}\Gamma(T^{*}M\otimes\Sigma^{g}M)\cong \Gamma(TM\otimes\Sigma^{g}M) \overset{\mu}{\to}\cal{F}^{g}.
	\]
 	The $\Spin_{n}$-representation $\bb{R}^{n}\otimes\Delta_{n}$ spilts as $\bb{R}^{n}\otimes\Delta_{n}=\ker(\mu)\oplus\Delta_{n}$ and this induces the decomposition $TM\otimes\Sigma^{g}M=\ker(\mu)\oplus \Sigma^{g}M$. We shall write  $p : \bb{R}^{n}\otimes\Delta_{n}\to\ker(\mu)$ for  the universal projection, which  locally is defined  by $ X\otimes\varphi\mapsto X\otimes\varphi +\frac{1}{n}\sum_{i}e_{i}\otimes e_{i}\cdot X\cdot\varphi$.   This naturally   extends to sections and yields the    {\it (generalized) Penrose, or twistor operator},  
	 \[
	 \cal{P}^{s}=p\circ\nabla^{s} :  \cal{F}^{g}\overset{\nabla^{s}}{\to}\Gamma(T^{*}M\otimes\Sigma^{g}M)\cong \Gamma(TM\otimes\Sigma^{g}M) \overset{p}{\to} \Gamma(\ker(\mu)).
	 \]
	%For $s=0$,  the  operators  $D^{0}\equiv D^{g}:=\mu\circ\nabla^{g}$ and $\cal{P}^{0}\equiv \cal{P}^{g}:=p\circ\nabla^{g}$ coincide with the  Riemannian Dirac operator  and the Riemannian   twistor operator, respectively. 
	 Locally, the operators $D^{s}$ and $\cal{P}^{s}$ attain the expressions 
	 \[
	 D^{s}(\varphi)=\sum_{i=1}^{n}e_{i}\cdot\nabla^{s}_{e_{i}}\varphi, \quad\text{and}\quad \cal{P}^{s}(\varphi):= \sum_{i=1}^{n} e_{i}\otimes\{\nabla^{s}_{e_{i}}\varphi+\frac{1}{n}e_{i}\cdot D^{s}(\varphi)\},
	 \]
	  respectively.   %=D^{g}(\varphi)+3sT\cdot\varphi=D^{g}(\varphi)+3s\sum_{i<j<k}T(e_{i}, e_{j}, e_{k})e_{i}\cdot e_{j}\cdot e_{k}\cdot\varphi,\\
   % \cal{P}^{s}(\varphi)&:=& \sum_{i=1}^{n} e_{i}\otimes\{\nabla^{s}_{e_{i}}\varphi+\frac{1}{n}e_{i}\cdot D^{s}\varphi\}=\cal{P}^{g}(\varphi)+s\sum_{i}e_{i}\otimes\Big\{\Big[(e_{i}\lrcorner T)+\frac{3}{n}e_{i}\cdot T\Big]\cdot\varphi\Big\}.
   % \end{eqnarray*}
       \noindent Finally we mention  that the one-parameter family of generalized Dirac operators $\{D^{s}\equiv D^{g}+3sT : s\in\bb{R}^{*}\}$ and the  Riemannian Dirac operator $D^{0}\equiv D^{g}$ are sharing several common properties.    For example,  $D^{s}$ if formally self-adjoint  
 in $L^{2}(\Sigma^{g}M)$ for any $s\in\bb{R}$, since the torsion $T^{s}=4sT$ is a 3-form  \cite{FS}. Moreover,  and similarly with  the Riemannian case $(s=0)$ (cf. \cite{Fr, Ginoux}), one can show that:
\begin{prop} \label{usef1}
 For any $s\in\bb{R}$, $f\in\cal{C}^{\infty}(M; \bb{R})$, $X\in\Gamma(TM)$, $\xi\in T^{*}M$, $\omega\in\Lambda^{p}T^{*}M$ and $\varphi\in\cal{F}^{g}$, the following hold: \\
 $(1)$  $D^{s}(f\varphi) = {\rm grad}(f) \cdot \varphi + f D^{s}(\varphi)$.\\
 $(2)$   The principal symbol of $D^{s}$ is given by $\sigma(D^{s})(\xi)(\varphi)=\xi^{\sharp}\cdot\varphi$ and hence $D^{s}$ is elliptic.\\
 $(3)$  The    operator $-(D^{s})^{2}$ is strongly elliptic, i.e. $\langle \sigma(-(D^{s})^{2})(\xi)\varphi, \varphi\rangle=|\xi|^{2}|\varphi|^{2}$.\\
  $(4)$ $ D^{s}(X\cdot \varphi)=\sum_{j}e_{j}\cdot (\nabla^{s}_{e_{j}}X)\cdot\varphi-X\cdot D^{s}(\varphi)-2\nabla^{s}_{X}\varphi$.\\
  $(5)$ $D^{s}(w\cdot\varphi)=(-1)^{p}\omega\cdot D^{s}(\varphi)+(d^{s}\omega+\delta^{s}\omega)\cdot\varphi-2\sum_{j}(e_{j}\lrcorner \omega)\cdot\nabla^{s}_{e_{j}}\varphi$.
   \end{prop}
\section{The generalized $\frac{1}{2}$-Ricci type  formula and basic applications}\label{1/2}
 
 \subsection{The generalized $\frac{1}{2}$-Ricci type  formula} In this section we shall introduce the generalized $\frac{1}{2}$-Ricci type formula. For this is useful to fix, once and for all, a  Riemannian spin manifold $(M^{n}, g, T)$ $(n\geq 3)$ endowed with a non-trivial 3-form $T\in\Lambda^{3}T^{*}M$, such that $\nabla^{c}T=0$, where $\nabla^{c}:=\nabla^{g}+\frac{1}{2}T$.    As we have already pointed out in Remark \ref{ricci},   the reason of our assumption $\nabla^{c}T=0$ is the symmetry  of the Ricci tensor $\Ric^{s}$ and Lemma \ref{JUL1}.  %\noindent We begin with the $\frac{1}{2}$-$\Ric^{s}$-identity. %As we explained in Introduction \ref{intro}, our approach to this result  is  different for example  from  \cite{FrIv}, or \cite{Agr03} and it relies on  the new   $\frac{1}{2}$-$\Ric^{s}$-formula, which we are going to prove below.      
%  \vskip 0.2cm
\begin{lemma} \label{LEMA}   \textnormal{\bf(The generalized $\frac{1}{2}$-Ricci type  formula, or $\frac{1}{2}$-$\Ric^{s}$-formula)}  Assume that $\nabla^{c}T=0$. Then,  the Ricci endomorphism  $\Ric^{s}(X)$  satisfies  
\begin{eqnarray*}
\frac{1}{2}\Ric^{s}(X)\cdot \varphi &=& D^{s}(\nabla^{s}_{X}\varphi)-\nabla^{s}_{X}(D^{s}\varphi)-\sum_{j=1}^{n}e_{j}\cdot \Big[\nabla^{s}_{\nabla^{s}_{e_{j}}X}\varphi+4s\nabla^{s}_{T(X, e_{j})}\varphi\Big]\nonumber\\
&&+s(3-4s)(X\lrcorner \sigma_{T})\cdot\varphi, 
\end{eqnarray*}
for any arbitrary  vector field $X\in\Gamma(TM)$,  spinor field $\varphi\in\cal{F}^{g}$  and $s\in\bb{R}$.
\end{lemma}

  \vskip 0.2cm
    \noindent {\bf Proof.} 
We use Lemma \ref{JUL1} and   replace $\cal{R}^{s}_{X, e_{j}}\varphi$ by $\nabla^{s}_{X}\nabla^{s}_{e_{j}}\varphi-\nabla^{s}_{e_{j}}\nabla^{s}_{X}\varphi-\nabla^{s}_{[X, e_{j}]}\varphi$. This yields the following
\begin{eqnarray*}
\frac{1}{2}\Ric^{s}(X)\cdot\varphi&=&-\sum_{j}e_{j}\cdot \cal{R}^{s}_{X, e_{j}}\varphi+s(3-4s)(X\lrcorner \sigma_{T})\cdot\varphi \\
&=&-\sum_{j}e_{j}\cdot \Big\{\nabla^{s}_{X}\nabla^{s}_{e_{j}}\varphi-\nabla^{s}_{e_{j}}\nabla^{s}_{X}\varphi-\nabla^{s}_{[X, e_{j}]}\varphi\Big\}+s(3-4s)(X\lrcorner \sigma_{T})\cdot\varphi \\
&=&\sum_{j}e_{j}\cdot (\nabla^{s}_{e_{j}}\nabla^{s}_{X}\varphi) -\sum_{j}e_{j}\cdot (\nabla^{s}_{X}\nabla^{s}_{e_{j}}\varphi)+\sum_{j}e_{j}\cdot \nabla^{s}_{[X, e_{j}]}\varphi+s(3-4s)(X\lrcorner \sigma_{T})\cdot\varphi,\\
&=&D^{s}(\nabla^{s}_{X}\varphi)-\sum_{j}\nabla^{s}_{X}(e_{j}\cdot\nabla^{s}_{e_{j}}\varphi)+\sum_{j}(\nabla^{s}_{X}e_{j})\cdot(\nabla^{s}_{e_{j}}\varphi)+\sum_{j}e_{j}\cdot\nabla^{s}_{[X, e_{j}]}\varphi\\
&&+s(3-4s)(X\lrcorner \sigma_{T})\cdot\varphi\\
&=&D^{s}(\nabla^{s}_{X}\varphi)-\nabla^{s}_{X}(D^{s}\varphi)+\sum_{j}(\nabla^{s}_{X}e_{j})\cdot(\nabla^{s}_{e_{j}}\varphi)+\sum_{j}e_{j}\cdot\nabla^{s}_{[X, e_{j}]}\varphi\\
&&+s(3-4s)(X\lrcorner \sigma_{T})\cdot\varphi,
\end{eqnarray*}
where for the fourth  equality we applied the Liebniz rule    to replace
$e_{j}\cdot(\nabla^{s}_{X}\nabla^{s}_{e_{j}}\varphi)$ by $\nabla^{s}_{X}(e_{j}\cdot \nabla^{s}_{e_{j}}\varphi)-(\nabla^{s}_{X}e_{j})\cdot(\nabla^{s}_{e_{j}}\varphi)$. By the definition of $T^{s}$, we also have that $[X, e_{j}]=\nabla^{s}_{X}e_{j}-\nabla^{s}_{e_{j}}X-T^{s}(X, e_{j})=\nabla^{s}_{X}e_{j}-\nabla^{s}_{e_{j}}X-4sT(X, e_{j})$ and consequently  
 $ \sum_{j}e_{j}\cdot\nabla^{s}_{[X, e_{j}]}\varphi=\sum_{j}e_{j}\cdot\Big[\nabla^{s}_{\nabla^{s}_{X}e_{j}}\varphi-\nabla^{s}_{\nabla^{s}_{e_{j}}X}\varphi-4s\nabla^{s}_{T(X, e_{j})}\varphi\Big]$. 
 Hence, the formula given above  reduces to the following one:
 \begin{eqnarray*}
\frac{1}{2}\Ric^{s}(X)\cdot\varphi&=&D^{s}(\nabla^{s}_{X}\varphi)-\nabla^{s}_{X}(D^{s}\varphi)-\sum_{j}e_{j}\cdot \nabla^{s}_{\nabla^{s}_{e_{j}}X}\varphi +\sum_{j}\Big[(\nabla^{s}_{X}e_{j})\cdot(\nabla^{s}_{e_{j}}\varphi)+e_{j}\cdot\nabla^{s}_{\nabla^{s}_{X}e_{j}}\varphi\Big]\\
&&-4s\sum_{j}e_{j}\cdot\nabla^{s}_{T(X, e_{j})}\varphi +s(3-4s)(X\lrcorner \sigma_{T})\cdot\varphi.
\end{eqnarray*}
Now,  it is easy to see that $\sum_{j}\Big[(\nabla^{s}_{X}e_{j})\cdot(\nabla^{s}_{e_{j}}\varphi)+e_{j}\cdot\nabla^{s}_{\nabla^{s}_{X}e_{j}}\varphi\Big]=0$,  for any $s\in\bb{R}$,  $X\in\Gamma(TM)$ and $\varphi\in\cal{F}^{g}$ (for the Riemannian case $s=0$, see also \cite[p.~133]{FKim}).  For example, this immediately follows after using (without loss of generality)   a local $\nabla^{s}$-parallel frame, i.e.   a local orthonormal frame  $\{e_j\}$ satisfying  $(\nabla^{s}e_{j})_{x}=0$ at $x\in M$, for any $j=1, \ldots, n$.  This finishes the proof. $\blacksquare$

\begin{remark}
\textnormal{Notice that the $\frac{1}{2}$-$\Ric^{s}$-formula (\ref{gen1/2ric})
can be  simplified a little further.  Indeed, in terms of our $\nabla^{s}$-parallel frame $\{e_{i}\}$ it  is $[X, e_{i}]=-\nabla^{s}_{e_{i}}X-4sT(X, e_{i})$ and the third term in (\ref{gen1/2ric}) reduces to
\[
-\sum_{j=1}^{n}e_{j}\cdot \Big[\nabla^{s}_{\nabla^{s}_{e_{j}}X}\varphi+4s\nabla^{s}_{T(X, e_{j})}\varphi\Big]=\sum_{j=1}^{n}e_{j}\cdot \nabla^{s}_{[X, e_{j}]}\varphi.
\]
Thus, an equivalent expression of the $\frac{1}{2}$-$\Ric^{s}$-identity  is given by 
  \begin{equation}
  \frac{1}{2}\Ric^{s}(X)\cdot \varphi = 
 D^{s}(\nabla^{s}_{X}\varphi)-\nabla^{s}_{X}(D^{s}\varphi)+\sum_{j=1}^{n}e_{j}\cdot  \nabla^{s}_{[X, e_{j}]}\varphi + s(3-4s)(X\lrcorner \sigma_{T})\cdot\varphi.\label{equiva} 
 \end{equation}
In introduction we chose to present   (\ref{gen1/2ric}), instead of (\ref{equiva}), since the reduction to (\ref{1/2FK}) for $s=0$ is direct. In fact, for the  particular scopes of this work, relation   (\ref{gen1/2ric})  turns out to be more flexible  than    (\ref{equiva}).}%However, later in some cases we shall use   (\ref{equiva}) to  simply the computations.}
\end{remark}

\subsection{The new proof of the generalized SL-formula.} 
\noindent 
 Let us  present now some basic applications of the new $\frac{1}{2}$-$\Ric^{s}$-formula, e.g.  an alternative  new proof of the following well-known and fundamental result  (here we state the version for $\nabla^{c}T=0$ only  and refer to  the cited articles for general skew-torsion). %  and   the interested reader can find the result for  general torsion to the cited references). 
 
 \begin{theorem}\label{parallel} \textnormal{(\cite[Thm.~3.1]{FrIv}, \cite[Thm.~6.1]{AF},  \cite[Thm.~2.1]{ABK})} \label{AGFR}
 Under the assumtpion $\nabla^{c}T=0$,  any spinor field $\varphi\in\cal{F}^{g}:=\Gamma(\Sigma^{g}M)$ on $(M^{n}, g, T)$ satisfies  the relation %the following holds (see  \cite[Thm.~6.1]{AF} or \cite[Thm.~2.1]{ABK})
\begin{equation}\label{GSL}
(D^{s})^{2}(\varphi)=\Delta^{s}(\varphi)+ s(3-4s) dT\cdot\varphi-4s\slashed{D}^{s}(\varphi)+\frac{1}{4}\Sca^{s}\cdot \ \varphi, 
\end{equation}
 where $\slashed{D}^{s}$ is the first-order differetntial operator   defined by $\slashed{D}^{s}:=\sum_{i}(e_{i}\lrcorner T)\cdot\nabla^{s}_{e_{i}}\varphi$ and 
 $ \Delta^{s}:=(\nabla^{s})^{*}\nabla^{s}:=-\sum_{i}\big[\nabla^{s}_{e_{i}}\nabla^{s}_{e_{i}}+\nabla^{s}_{\nabla^{g}_{e_{i}}e_{i}}\big]$  denotes the spin Laplace operator associated to $\nabla^{s}$.
 \end{theorem}
\noindent First we recall  that
\begin{lemma}\label{USEF3} Consider a  $p$-form $\omega\in\Lambda^{p}T^{*}M$ and an orthonormal frame $\{e_{j}\}$.  Then  
\[
\sum_{j=1}^{n}e_{j}\cdot(e_{j}\lrcorner\omega)=p\omega, \quad \sum_{j=1}^{n}e_{j}\cdot (e_{j}\wedge \omega)=(p-n)\omega.
\]
\end{lemma}

\vskip 0.2cm
 \noindent 
{\bf New  proof of Theorem \ref{parallel}.}  Since the Ricci tensor $\Ric^{s}$ is symmetric, as in the Riemannian case,   one can use the generalized $\frac{1}{2}$-Ricci formula  and apply   a contraction with respect to a $\nabla^{s}$-parallel local orthonormal frame $\{e_{i}\}$. This means 
\[
\sum_{i}e_{i}\cdot \Ric^{s}(e_{i})\cdot\varphi=\sum_{i, k}\Ric^{s}(e_{i}, e_{k})\cdot e_{i}\cdot e_{k}\cdot\varphi=-\sum_{i}\Ric^{s}(e_{i}, e_{i})\cdot\varphi=-\Sca^{s}\cdot \ \varphi
\]
and hence, for the left-hand side part of (\ref{gen1/2ric}) we obtain $\frac{1}{2}\sum_{i}e_{i}\cdot \Ric^{s}(e_{i})\cdot\varphi=-\frac{1}{2}\Sca^{s}\cdot \ \varphi$. Next we focus   on the right-hand side and   write all together: 
\begin{eqnarray*}
-\frac{1}{2}\Sca^{s}\cdot \ \varphi&=&\sum_{i}e_{i}\cdot\Big[D^{s}(\nabla^{s}_{e_{i}}\varphi)-\nabla^{s}_{e_{i}}(D^{s}\varphi)\Big]  -\sum_{i, j}e_{i}\cdot e_{j}\cdot  \Big[\nabla^{s}_{\nabla^{s}_{e_{j}}e_{i}}\varphi+4s\nabla^{s}_{T(e_{i}, e_{j})}\varphi\Big]\\
&& +s(3-4s)\sum_{i} e_{i}\cdot (e_{i}\lrcorner \sigma_{T})\cdot\varphi\\
&=&\sum_{i}e_{i}\cdot D^{s}(\nabla^{s}_{e_{i}}\varphi)-\sum_{i}e_{i}\cdot \nabla^{s}_{e_{i}}(D^{s}\varphi)-\sum_{i, j}e_{i}\cdot e_{j}\cdot  \Big[\nabla^{s}_{\nabla^{s}_{e_{j}}e_{i}}\varphi+4s\nabla^{s}_{T(e_{i}, e_{j})}\varphi\Big]\\
&&+4s(3-4s)\sigma_{T}\cdot\varphi,\\
&=&\sum_{i}e_{i}\cdot D^{s}(\nabla^{s}_{e_{i}}\varphi)-(D^{s})^{2}(\varphi)-\sum_{i, j}e_{i}\cdot e_{j}\cdot\Big[\nabla^{s}_{\nabla^{s}_{e_{j}}e_{i}}\varphi+4s\nabla^{s}_{T(e_{i}, e_{j})}\varphi\Big]\\
&&+4s(3-4s)\sigma_{T}\cdot\varphi,
\end{eqnarray*}
where  we used the fact that  $\sum_{i}e_{i}\cdot (e_{i}\lrcorner \sigma_{T})=4\sigma_{T}$ (see Lemma \ref{USEF3}) and $-\sum_{i}e_{i}\cdot \nabla^{s}_{e_{i}}(D^{s}\varphi)=-(D^{s})^{2}(\varphi)$.  We proceed   with the   two sums appearing in the resulting formula.
For the first one, we are based on  the formula given in Proposition \ref{usef1}, $(4)$; we  replace $X$ by $e_{i}$ and $\varphi$ by $\nabla^{s}_{e_{i}}\varphi$ and this yields
\[
D^{s}(e_{i}\cdot \nabla^{s}_{e_{i}}\varphi)=\sum_{i} e_{i}\cdot (\nabla^{s}_{e_{i}}e_{i})\cdot\nabla^{s}_{e_{i}}\varphi-e_{i}\cdot D^{s}(\nabla^{s}_{e_{i}}\varphi)-2\nabla^{s}_{e_{i}}\nabla^{s}_{e_{i}}\varphi=
-e_{i}\cdot D^{s}(\nabla^{s}_{e_{i}}\varphi)-2\nabla^{s}_{e_{i}}\nabla^{s}_{e_{i}}\varphi.
\]
Consequently,
\[
\sum_{i}e_{i}\cdot D^{s}(\nabla^{s}_{e_{i}}\varphi)=-\sum_{i}D^{s}(e_{i}\cdot \nabla^{s}_{e_{i}}\varphi)-2\sum_{i}\nabla^{s}_{e_{i}}\nabla^{s}_{e_{i}}\varphi=-(D^{s})^{2}(\varphi)+2\Delta^{s}(\varphi),
\]
where  in terms of the $\nabla^{s}$-parallel local orthonormal frame $\{e_{i}\}$, the spin Laplace operator  reads by $\Delta^{s}(\varphi)=(\nabla^{s})^{*}\nabla^{s}\varphi=-\sum_{i}\nabla^{s}_{e_{i}}\nabla^{s}_{e_{i}}\varphi$ (observe  that the relation $(\nabla^{s}e_{i})_{x}=0$, yields $(\nabla^{g}_{e_{i}}e_{i})_{x}=0$).  Therefore, in order to complete the proof of Theorem \ref{parallel}, we  just need to show that  the  second  sum induces  the   operator $\slashed{D}^{s}$  with the desired coefficient.  Indeed, because  $(\nabla^{s}_{e_{j}}e_{i})_{x}=0$, we obtain  
\[
-\sum_{i, j}e_{i}\cdot e_{j}\cdot\Big[\nabla^{s}_{\nabla^{s}_{e_{j}}e_{i}}\varphi+4s\nabla^{s}_{T(e_{i}, e_{j})}\varphi\Big]=-4s\sum_{i, j}e_{i}\cdot e_{j}\cdot\nabla^{s}_{T(e_{i}, e_{j})}\varphi.
\]
%where of course  the torsion $T$ is not trivial even  locally. %  in terms of the $\nabla^{s}$-parallel local orthonormal frame, it is $T(e_{i}, e_{j})=-\frac{1}{4s}[e_{i}, e_{j}]\neq 0$.
   For a few, we forgot the factor  $-4s$ and since $T(e_i, e_j)=\sum_{k}T^{k}_{ij}e_{k}=\sum_{k}T(e_{i}, e_{j}, e_{k})e_{k}$, it follows that 
   \begin{eqnarray}
\sum_{i, j}e_{i}\cdot e_{j}\cdot\nabla^{s}_{T(e_{i}, e_{j})}\varphi&=&\sum_{i, j}e_{i}\cdot e_{j}\cdot\nabla^{s}_{\sum_{k}T_{ij}^{k}e_{k}}\varphi\nonumber\\
&=& \sum_{i, j, k}T^{k}_{ij}e_{i}\cdot e_{j}\cdot\nabla^{s}_{e_{k}}\varphi=\sum_{i, j, k}T(e_{i}, e_{j}, e_{k})e_{i}\cdot e_{j}\cdot\nabla^{s}_{e_{k}}\varphi, \label{classic}
   \end{eqnarray}
   for some non-zero real numbers $T^{k}_{ij}:=T(e_{i}, e_{j}, e_{k})=g\big(T(e_{i}, e_{j}), e_{k}\big)$, with $T^{k}_{ij}=-T^{k}_{ji}=-T^{i}_{jk}$, etc. Recall now that    $\sum_{j}T(X, e_{j})\cdot e_{j}=-2(X\lrcorner T)$ for any $X\in \Gamma(TM)$ (see for example \cite[p.~328]{ABK}). Thus,  
\begin{eqnarray}
-2(e_{k}\lrcorner T)&=&\sum_{j}T(e_{k}, e_{j})\cdot e_{j}=\sum_{i, j}g\big(T(e_{k}, e_{j}), e_{i}\big) e_{i}\cdot e_{j}\nonumber\\
&=&\sum_{i, j}T(e_{k}, e_{j}, e_{i})e_{i}\cdot e_{j}=-\sum_{i, j}T(e_{i}, e_{j}, e_{k})e_{i}\cdot e_{j},\label{ekT}
\end{eqnarray}
and  a combination with (\ref{classic}) yields  
\begin{equation}\label{motiv}
\sum_{i, j}e_{i}\cdot e_{j}\cdot\nabla^{s}_{T(e_{i}, e_{j})}\varphi=2\sum_{k}(e_{k}\lrcorner T)\cdot \nabla^{s}_{e_{k}}\varphi=2\slashed{D}^{s}(\varphi).
\end{equation}
   Adding all our results together we obtain  
  \begin{eqnarray*}
  -\frac{1}{2}\Sca^{s}\cdot \ \varphi=-2(D^{s})^{2}(\varphi)+2\Delta^{s}(\varphi)-8s\slashed{D}^{s}(\varphi)+4s(3-4s)\sigma_{T}\cdot\varphi.
  \end{eqnarray*}
  Since $\nabla^{c}T=0$ it is $dT=2\sigma_{T}$ and consequently  the last identity is nothing than the generalized Schr\"odinger-Lichnerowicz formula under the condition $\nabla^{c}T=0$. $\blacksquare$

       \subsection{The new proof of the twistorial $\frac{1}{2}$-$\Ric^{s}$-formula}\label{tsts}
      \noindent  Here we shall  describe the restriction of the    $\frac{1}{2}$-$\Ric^{s}$-formula to twistor spinors, providing a simpler  proof of the twistorial $\frac{1}{2}$-$\Ric^{s}$-formula. This has been recently introduced  by the author in \cite[Lem.~2.2]{Chrysk2}, with different however methods.   Recall that a {\it twistor spinor with torsion} (TsT in short)   is a   (non-trivial) spinor field $\varphi\in\cal{F}^{g}$, solving  the  equation 
  \[
  \nabla^{s}_{X}\varphi+\frac{1}{n}X\cdot D^{s}(\varphi)=0, \quad \text{for  any} \   X\in\Gamma(TM),
  \]
  for some parameter $s\neq 0$. Hence, a TsT  is a  spinor field belonging to  the kernel of the   generalized  twistor operator  $\cal{P}^{s}$,  see \cite{ABK, Chrysk2} for more details.  Obviously, for $s=0$ one obtains the usual notion of Riemannian twistor spinors, i.e. elements  $\varphi\in\ker\cal{P}^{g}$, where $\cal{P}^{g}\equiv\cal{P}^{0}$ (cf. \cite{Fr, Ginoux}).
  \begin{theorem}\label{prof}
 \textnormal{\bf (Twistorial $\frac{1}{2}$-$\Ric^{s}$-formula)}
 The  action of the $\frac{1}{2}$-$\Ric^{s}$-formula on a non-trivial twistor spinor   $\varphi\in\ker(\cal{P}^{s})$ is given by   
 \[
 \frac{1}{2}\Ric^{s}(X)\cdot\varphi=\frac{1}{n}X\cdot (D^{s})^{2}(\varphi) -\frac{n-2}{n}\nabla^{s}_{X}\big(D^{s}(\varphi)\big)+\frac{8s}{n}(X\lrcorner T)\cdot D^{s}(\varphi) +s(3-4s)(X\lrcorner \sigma_{T})\cdot\varphi,
 \]
  for any $X\in\Gamma(TM)$.  In particular, for $s=0$ and  for a Riemannian twistor spinor, it induces the relation 
 \[
 \nabla^{g}_{X}\big(D^{g}(\varphi)\big)=\frac{n}{2(n-2)}\Big[-\Ric^{g}(X)\cdot\varphi+\frac{\Sca^{g}}{2(n-1)}X\cdot\varphi\Big]=\frac{n}{2}\Scho^{g}(X)\cdot\varphi,\]
   where $\Scho^{g}(X):=\frac{1}{n-2}\big[-\Ric^{g}(X)\cdot\varphi+\frac{\Sca^{g}}{2(n-1)}X\big]$ is the  Schouten endomorphism  associated to  $\nabla^{g}$.
 \end{theorem}

%\begin{theorem}
%%The restriction of the $\frac{1}{2}$-$\Ric^{s}$-formula (\ref{gen1/2ric})   to the kernel   $\ker(\cal{P}^{s})\subset\cal{F}^{g}$ for some  parameter $s$, yields the Proposition \ref{good}, i.e. the twistorial $\frac{1}{2}$-$\Ric^{s}$-identity. 
%\end{theorem} 
\vskip 0.2cm
\noindent {\bf Proof.} 
 We apply the $\frac{1}{2}$-$\Ric^{s}$-formula (\ref{gen1/2ric}) to a non-trivial twistor spinor $\varphi\in\ker(\cal{P}^{s})$. Any vector field $X\in\Gamma(TM)$ satisfies  $\nabla^{s}_{X}\varphi=-\frac{1}{n}X\cdot D^{s}(\varphi)$, hence for the first term in (\ref{gen1/2ric}), Proposition \ref{usef1}, (4), yields that 
 \begin{eqnarray*}
 D^{s}(\nabla^{s}_{X}\varphi)&=&D^{s}(-\frac{1}{n}X\cdot D^{s}(\varphi))=-\frac{1}{n}D^{s}(X\cdot D^{s}(\varphi))\\
 &=&-\frac{1}{n}\Big[\sum_{j}e_{j}\cdot (\nabla^{s}_{e_{j}}X)\cdot D^{s}(\varphi)-X\cdot D^{s}(D^{s}(\varphi))-2\nabla^{s}_{X}(D^{s}(\varphi))\Big]\\
 &=&-\frac{1}{n}\sum_{j}e_{j}\cdot (\nabla^{s}_{e_{j}}X)\cdot D^{s}(\varphi)+\frac{1}{n}X\cdot (D^{s})^{2}(\varphi)+\frac{2}{n}\nabla^{s}_{X}(D^{s}(\varphi)).
 \end{eqnarray*}
 Therefore,   for the first two terms of (\ref{gen1/2ric}) we deduce that
   \begin{equation}\label{perf}
D^{s}(\nabla^{s}_{X}\varphi)-\nabla^{s}_{X}(D^{s}(\varphi))=\frac{1}{n}X\cdot (D^{s})^{2}(\varphi)-\frac{n-2}{n}\nabla^{s}_{X}\big(D^{s}(\varphi)\big)
-\frac{1}{n}\sum_{j}e_{j}\cdot (\nabla^{s}_{e_{j}}X)\cdot D^{s}(\varphi).  %\quad (\AcPa )
 \end{equation}
So, we  have already constructed  two desired terms with the right coefficients. Now, let us consider   the sum $
-\sum_{j=1}^{n}e_{j}\cdot \Big[\nabla^{s}_{\nabla^{s}_{e_{j}}X}\varphi+4s\nabla^{s}_{T(X, e_{j})}\varphi\Big]$  for some  non-trivial $\varphi\in\ker(\cal{P}^{s})$.    Since $\nabla^{s}_{\nabla^{s}_{e_{j}}X}\varphi=-\frac{1}{n}(\nabla^{s}_{e_{j}}X)\cdot D^{s}(\varphi)$, it follows that 
\[
-\sum_{j=1}^{n}e_{j}\cdot \nabla^{s}_{\nabla^{s}_{e_{j}}X}\varphi=\frac{1}{n}\sum_{j=1}^{n}e_{j}\cdot (\nabla^{s}_{e_{j}}X)\cdot D^{s}(\varphi), 
\]
and this is canceled with the third term in the right-hand side of $(\ref{perf})$. Moreover,  based on the identity $2(X\lrcorner T)=\sum e_{j}\cdot T(X, e_{j})$ we obtain
\[
-4s\sum_{j=1}^{n}e_{j}\cdot \nabla^{s}_{T(X, e_{j})}\varphi=\frac{4s}{n}\sum_{j=1}^{n}e_{j}\cdot T(X, e_{j})\cdot D^{s}(\varphi) 
=\frac{8s}{n}(X\lrcorner T)\cdot D^{s}(\varphi)
\]
and our claim follows. For more details related to the case $s=0$ we refer to \cite{Fr, Ginoux, Chrysk2}. $\blacksquare$
\vskip 0.2cm
\noindent The twistorial $\frac{1}{2}$-$\Ric^{s}$-formula  gives rise to
\[
 \frac{1}{2}\Sca^{s}\cdot \ \varphi = -\frac{24s}{n}T\cdot D^{s}(\varphi)+\frac{2(n-1)}{n}(D^{s})^{2}(\varphi)-4s(3-4s)\sigma_{T}\cdot\varphi,
\]
which it is not hard to see that is  equivalent with the generalized SL-formula associated to the   Dirac operator $D^{s}$, when  this operator is  restricted to  $\ker(\cal{P}^{s})$.  The most basic  consequences of the identity stated  in Theorem \ref{prof}, have been described  in \cite{Chrysk2}. For example:
\begin{prop} \label{lab1} \textnormal{(\cite{Chrysk2})}
Let $(M^{n}, g, T)$ $(n\geq 3)$ be a connected Riemannian spin manifold with $\nabla^{c}T=0$. Then,

\noindent   a) The kernel of the twistor operator $\cal{P}^{s}$ is   finite dimensional, i.e. $\dim_{\bb{C}}\ke(\cal{P}^{s})\leq 2^{[\frac{n}{2}]+1}$.
 
\noindent b) If  $\varphi$ and $D^{s}(\varphi)$ vanish at some point $p\in M$ and $\varphi\in\ke(\cal{P}^{s})$, then   $\varphi\equiv 0$.
\end{prop}

 \noindent   Obviously, Proposition \ref{lab1} generalizes  classical properties  of  Riemannian twistor spinors (cf. \cite{Fr, Ginoux}), to the whole family $\{\nabla^{s} : s\in\bb{R}\}$.  Notice however that in the Riemannian case   the space $\ker(\cal{P}^{g})$ is  in addition  a  conformal invariant of $(M^{n}, g)$. 
 A similar result for connections with skew-torsion     is  known to hold only in dimension 4 (cf.  \cite{Dalakov}). For more details on TsT and interesting examples of special geometric structures admitting  this kind of spinor fields, we refer to  \cite{Julia, ABK, Chrysk2} (see also Section \ref{pKsT} below).

 \section{$\nabla^{s}$-parallel spinors and $\nabla^{c}$-parallel spinors}\label{psps}
 \subsection{Parallel spinors}
 \noindent  In this section we present    applications of the  $\frac{1}{2}$-$\Ric^{s}$-formula,  related with $\nabla^{s}$-parallel spinors. Since Lemma \ref{LEMA} holds only under the assumption $\nabla^{c}T=0$, it should be pointed out (even if this is not repeated throughout), that  these results are meant for spin manifolds $(M^{n}, g, T)$ with $\nabla^{c}T=0$.  Hence, fix a Riemannian spin manifold $(M^{n}, g)$ endowed with a non-trivial 3-form $T\in\Lambda^{3}T^{*}M$ such that $\nabla^{c}T=0$, and consider the lift of the 1-parameter family $\{\nabla^{s}=\nabla^{g}+2sT : s\in\bb{R}\}$ to  the spinor bundle $\Sigma^{g}M$.  
  By a {\it $\nabla^{s}$-parallel spinor} we mean a non-trivial spinor field $\varphi_{0}\in\cal{F}^{g}$  satisfying the equation $\nabla^{s}_{X}\varphi_{0}=0$, for some $s\in\bb{R}$ and any $X\in\Gamma(TM)$. This notion includes the following two well-known kinds of parallel spinors:
   
 \begin{itemize}
 \item   $s=0$;  then we speak about $\nabla^{g}$-parallel spinors and their existence yields the $\Ric^{g}$-flatness of $(M^{n}, g)$, i.e. $\Ric^{g}(X)\cdot\varphi_{0}=0$ for any $X\in\Gamma(TM)$, see for example \cite{Fr, Ginoux} (notice that an easy way to prove the Ricci flatness  is via the $\frac{1}{2}$-$\Ric^{g}$-formula (\ref{1/2FK})).
  \item  $s=1/4$;  then we speak about $\nabla^{c}$-parallel spinors and is a simple consequence of \cite[Col.~3.2]{FrIv} that when the condition $\nabla^{c}T=0$ holds, then a solution of the relation $\nabla^{c}\varphi_{0}=0$ must satisfies the relations $\Ric^{c}(X)\cdot\varphi_{0}=(X\lrcorner \sigma_{T})\cdot\varphi_{0}$, for any $X\in\Gamma(TM)$ and $\Sca^{c}\cdot \ \varphi_{0}=-4\sigma_{T}\cdot\varphi_{0}$.
 \end{itemize}

\vskip 0.1cm
\noindent  The $\frac{1}{2}$-$\Ric^{s}$-formula immediately  yields  integrability conditions for any member of the family $\{\nabla^{s} : s\in\bb{R}\}$.  Moreover, when a $\nabla^{c}$-parallel spinor exists   it allows us  to describe the action of the endomorphism $\Ric^{s}(X) : \cal{F}^{g}\to\cal{F}^{g}$ for any other $s$.   We begin with the following:
  \begin{corol}\label{genpar}
  Assume that $\nabla^{c}T=0$ and let $\varphi_{0}\in\cal{F}^{g}$ be a non-trivial spinor field which is parallel with respect to $\nabla^{s}$, for some $s\in\bb{R}$. Then,   for the same $s$ and for any $X\in\Gamma(TM)$ the spinor  $\varphi_{0}$ must satisfy the following:
  \begin{eqnarray}
  \Ric^{s}(X)\cdot\varphi_{0}&=&2s(3-4s)(X\lrcorner \sigma_{T})\cdot\varphi_{0},\label{sparallel}\\
    \Sca^{s}\cdot \ \varphi_{0}&=&-8s(3-4s)\sigma_{T}\cdot\varphi_{0}.\label{scall2}
  \end{eqnarray}
\end{corol}

\begin{remark}\textnormal{Notice that the connection $\nabla^{3/4}$ has torsion $3T$ and (\ref{sparallel}) shows that  the existence of a $\nabla^{3/4}$-parallel spinor $\varphi_{0}$ implies the $\Ric^{3/4}$-flatness of  $(M^{n}, g, T)$.}
\end{remark}
\noindent   For the record, we also mention that
\begin{corol}\label{riccflat} 
Assume that $\nabla^{c}T=0$  and let    $\varphi_{0}$ be a non-trivial  $\nabla^{s}$-parallel spinor for some  $s\in \bb{R}\backslash\{0, 3/4\}$. Then, $(M, g, T)$   is $\Ric^{s}$-flat for the same parameter $s$, if and only if $(X\lrcorner \sigma_{T})\cdot\varphi_{0}=0$ for any $X\in\Gamma(TM)$. 
 \end{corol}
 
  \vskip 0.2cm
 \noindent Finally, since the 4-form $\sigma_{T}$ vanishes in any dimension $n\leq 4$ (\cite{FrIv, AF}), we   have that
 \begin{corol}\label{rfl}
 A 3-dimensional or 4-dimensional Riemannian spin manifold $(M^{n}, g, T)$ with $\nabla^{c}T=0$, which admits a non-trivial $\nabla^{s}$-parallel spinor  $\varphi_{0}\in\ker(\nabla^{s})$ for some $s\in\bb{R}$,  is $\Ric^{s}$-flat for the same parameter $s$.
 \end{corol}
 \begin{example}
    \textnormal{Consider  the round 3-sphere $(\Ss^{3}, g_{\rm can})$ endowed with the volume form $T={\rm Vol}_{\Ss^{3}}$. The characteristic connection $\nabla^{c}=\nabla^{\pm 1/4}$  (which is not unique because $\Ss^{3}\cong\Spin_{3}\cong\SU_{2}$ is a Lie group), is induced by the Killing spinor equation.  The real   Killing spinors  of $\Ss^{3}$ trivialize its spinor bundle and  they are  $\nabla^{c}$-parallel, see  \cite[p.~729]{AF}. Hence,  Corollary \ref{genpar} or Corollary \ref{rfl} apply, and show  that any such spinor $\{\varphi_{j} :  1\leq j\leq 2^{[\frac{3}{2}]}\}$ must satisfy the equation  $\Ric^{c}(X)\cdot\varphi_{j}=0$ for any $X\in\Gamma(T\Ss^{3})$, for another approach see for example \cite[Prop.~5.1 and p.~133]{Chrysk2}.  More general, any simply connected compact Lie group $G$ with a bi-invariant metric $g$ is flat with respect to the Cartan-Schouten connections $\nabla^{\pm 1/4}$ and there are $\nabla^{\pm 1/4}$-parallel spinors which  satisfy $\Ric^{\pm 1/4}(X)\cdot\varphi=0$ for any $X\in\fr{g}$. Further   $\Ric^{c}$-flat structures carrying $\nabla^{c}$-parallel spinors can be found in \cite{AF, FrIv2,  FriedG2f}, for instance. }
    \end{example}

\begin{remark}\textnormal{(\cite{AF})} \textnormal{In the compact case, Agricola and Friedrich  \cite[Thm.~7.1]{AF}  proved  that there are at most three parameters with $\nabla^{s}$-parallel spinors.   Indeed, assume that $(M^{n}, g, T)$  is  a compact  Riemannian spin manifold endowed with a non-trivial 3-form $T\in\Lambda^{3}T^{*}M$, such that $\nabla^{c}T=0$. Then, any $\nabla^{s}$-parallel spinor $\varphi_{0}$ of unit length satisfies
 \[
8s(3-4s) \int_{M}\langle\sigma_{T}\cdot\varphi_{0}, \varphi_{0}\rangle  v^{g}+\int_{M}\Sca^{s} v^{g}=0.
 \]
This follows immediately  by integrating the   condition  (\ref{scall2}). Based now on (\ref{clasdif1}) one can show that if the mean value of $\langle\sigma_{T}\cdot\varphi_{0}, \varphi_{0}\rangle$ does not vanish, then the parameter $s$ equals to $s=1/4$, i.e. $\varphi_{0}$ is necessarily parallel under the characteristic connection. If the mean value of $\langle\sigma_{T}\cdot\varphi_{0}, \varphi_{0}\rangle$ vanishes, then the parameter $s$ depends on $\Sca^{g}$ and $\|T\|^{2}$. We refer to \cite{AF} for further details and examples. }
   \end{remark}
   %  \vskip 0.1cm
          \noindent  In the following we shall use  the $\frac{1}{2}$-$\Ric^{s}$-formula to describe the spinorial action of the Ricci endomorphism $\Ric^{s}(X)$ for any $s\in\bb{R}$, when there exists some $\nabla^{c}$-parallel spinor $\varphi_{0}$, without assuming however that $\nabla^{c}$ is the characteristic connection of some underlying special structure.          An important fact for our approach is that  the  torsion $T$ can be viewed as a ($\nabla^{c}$-parallel) symmetric endomorphism  on $\Sigma^{g}M$  in the sense that  \begin{equation}\label{symT}
\langle T\cdot 
 \varphi, \psi\rangle=\langle \varphi, T\cdot\psi\rangle, \quad \forall \ \varphi, \psi\in\cal{F}^{g}.
 \end{equation}
 Hence it is diagonalizable  with real eigenvalues.   % (functions).
 Then, one may decompose the spinor bundle  $\Sigma^{g}M$ into a direct sum of $T$-eigenbundles preserved by $\nabla^{c}$,  i.e. $\Sigma^{g}M=\bigoplus_{\gamma\in{\rm Spec}(T)}\Sigma^{g}_{\gamma}M$  with $\nabla^{c}\Sigma^{g}_{\gamma}M\subset \Sigma^{g}_{\gamma}M$. This induces a splitting also to  the space of sections, $ \cal{F}^{g}=\bigoplus_{\gamma\in{\rm Spec}(T)}\cal{F}^{g}(\gamma)$ with $\cal{F}^{g}(\gamma):=\Gamma(\Sigma^{g}_{\gamma}M)$. 
 We finally  remind that  when the torsion  $T$ is $\nabla^{c}$-parallel,  then any non-trivial $\nabla^{c}$-parallel spinor field  has constant $T$-eigenvalues, i.e. the equations  $\nabla^{c}T=0, \ \nabla^{c}\varphi_{0}=0$ and $T\cdot \varphi_{0}=\gamma\cdot \varphi_{0}$ for some $\gamma\in\Spec(T)$ imply that    $\gamma=\text{constant}\in\bb{R}$, see \cite[Thm.~1.1]{AFNP}.     
      
          \begin{theorem}\label{newpar}
 Consider   a Riemannian spin manifold $(M^{n}, g, T)$ $(n\geq 3)$  endowed with a non-trivial 3-form $T\in\Lambda^{3}T^{*}M$, such that   $\nabla^{c}T=0$, where $\nabla^{c}:=\nabla^{g}+\frac{1}{2}T$.  Assume that  $\varphi_{0}$ is a non-trivial $\nabla^{c}$-parallel spinor field  lying in $\cal{F}^{g}(\gamma)$, for some (constant) $\gamma\in\bb{R}$.  Then, for {\it any} $s\in\bb{R}$ and $X\in\Gamma(TM)$ the following holds
    \begin{eqnarray}
    \Ric^{s}(X)\cdot\varphi_{0}&=& -\frac{(16s^{2}-1)}{4}\sum_{j}e_{j}\cdot \big(T(X, e_{j})\lrcorner T\big)\cdot\varphi_{0} +\frac{(16s^2 + 3)}{4}(X\lrcorner \sigma_{T})\cdot\varphi_{0}\nonumber\\
    &=&\frac{(16s^{2}-1)}{4} \sum_{j} T(X, e_{j})\cdot (e_{j}\lrcorner T)\cdot\varphi_{0}+\frac{(16s^2 + 3)}{4}(X\lrcorner \sigma_{T})\cdot\varphi_{0}.\label{NPAR}
        \end{eqnarray}
             \end{theorem}

                          \vskip 0.2cm
          \noindent {\bf Proof}. The proof is rather long and relies on the $\frac{1}{2}$-$\Ric^{s}$-formula (\ref{gen1/2ric}) and the $\nabla^{c}$-parallelism of $T$.   To begin with,  notice  that a $\nabla^{c}$-parallel spinor $\varphi_{0}$ satisfies the following two equations (cf. \cite{Chrysk2})  
          \begin{equation}\label{usef3}
\nabla^{s}_{X}\varphi_{0}=\frac{4s-1}{4}(X\lrcorner T)\cdot\varphi_{0}, \quad\text{and}\quad D^{s}(\varphi_{0})=\frac{3(4s-1)}{4} \ T\cdot\varphi_{0},
\end{equation}
for any $X\in\Gamma(TM)$. In particular $\varphi_{0}$   is a $D^{s}$-eigenspinor, $D^{s}(\varphi_{0})=\frac{3(4s-1)\gamma}{4}\varphi_{0}$, for any $s\in\bb{R}$, with  $D^{c}(\varphi_{0})=0$ of course.  Let us  apply the $\frac{1}{2}$-$\Ric^{s}$-formula to $\varphi_{0}$. Due to (\ref{usef3}) and Proposition \ref{usef1}, (5),  for the first  term  of  (\ref{gen1/2ric}) we deduce that
\begin{eqnarray}
D^{s}(\nabla^{s}_{X}\varphi_{0})&=&\frac{(4s-1)}{4}D^{s}\big((X\lrcorner T)\cdot\varphi_{0}\big)\nonumber\\
&=&\frac{(4s-1)}{4}\Big[(X\lrcorner T)\cdot D^{s}(\varphi_{0})+\big(d^{s}+\delta^{s}\big)(X\lrcorner T)\cdot\varphi_{0}-2\sum_{j}(e_{j}\lrcorner X\lrcorner T)\cdot\nabla^{s}_{e_{j}}\varphi_{0}\Big]\nonumber\\
&=&3\Big[\frac{(4s-1)}{4}\Big]^{2}(X\lrcorner T)\cdot T\cdot\varphi_{0}+\frac{(4s-1)}{4}\Big[d^{s}(X\lrcorner T)+\delta^{s}(X\lrcorner T)\Big]\cdot\varphi_{0}\nonumber\\
&&-2\Big[\frac{(4s-1)}{4}\Big]^{2}\sum_{j}T(X, e_{j})\cdot (e_{j}\lrcorner T)\cdot\varphi_{0},\nonumber
%&=&3\Big[\frac{(4s-1)}{4}\Big]^{2}(X\lrcorner T)\cdot T\cdot\varphi_{0}+\frac{(4s-1)}{4}\Big[d^{s}(X\lrcorner T)+\delta^{s}(X\lrcorner T)\Big]\cdot\varphi_{0}\nonumber\\
%&&-6\Big[\frac{(4s-1)}{4}\Big]^{2}(X\lrcorner T)\cdot T\cdot\varphi_{0}\nonumber\\
%&=&-3\Big[\frac{(4s-1)}{4}\Big]^{2}(X\lrcorner T)\cdot T\cdot\varphi_{0}+\frac{(4s-1)}{4}d^{s}(X\lrcorner T)\cdot\varphi_{0}+\frac{(4s-1)}{4}\delta^{s}(X\lrcorner T)\cdot\varphi_{0}.\nonumber\\
%&=&\Big[\frac{(4s-1)}{4}\Big]^{2}\sum_{j}e_{j}\cdot(X\lrcorner T)\cdot (e_{j}\lrcorner T)\cdot\varphi_{0}+\frac{(4s-1)}{4}\Big[d^{s}(X\lrcorner T)+\delta^{s}(X\lrcorner T)\Big]\cdot\varphi_{0}.\nonumber
%&&\label{form1} 
%&&\label{form1} 
\end{eqnarray}
   where one  identifies  the 1-form $T(X, e_{j})^{\flat}:=g(T(X, e_{j}), -)=T(X, e_{j}, -)=e_{j}\lrcorner X\lrcorner T\in \Lambda^{1}T^{*}M$ with its dual vector field $T(X, e_{j})\in\Gamma(TM)$, via the metric tensor $g$.    From  (\ref{usef3}) and since $\gamma$ is a constant, we also obtain 
   \[
   \nabla^{s}_{X}\big(D^{s}(\varphi_{0})\big)=3\Big[\frac{(4s-1)}{4}\Big]^{2}\gamma\cdot(X\lrcorner T)\cdot\varphi_{0}=3\Big[\frac{(4s-1)}{4}\Big]^{2}  (X\lrcorner T)\cdot T\cdot \varphi_{0}
   \]   %Recall now that $-2(X\lrcorner \sigma_{T})=   (X\lrcorner T)\cdot T-T\cdot (X\lrcorner T)$, see \cite[p.~325]{ABK}, or \cite[p.~129]{Chrysk2}. 
  and  a combination with  the stated  expression for $D^{s}(\nabla^{s}\varphi_{0})$,  yields the difference %$\cal{A}(\varphi_{0}):=D^{s}(\nabla^{s}_{X}\varphi_{0})-\nabla^{s}_{X}\big(D^{s}(\varphi_{0})\big)$, i.e.
\begin{equation}\label{form3}
 D^{s}(\nabla^{s}_{X}\varphi_{0})-\nabla^{s}_{X}\big(D^{s}(\varphi_{0})\big)=-\frac{(4s-1)^{2}}{8}\sum_{j}T(X, e_{j})\cdot (e_{j}\lrcorner T)\cdot\varphi_{0}+\frac{(4s-1)}{4}\Big[d^{s}(X\lrcorner T)+\delta^{s}(X\lrcorner T)\Big]\cdot\varphi_{0}.
  \end{equation}
   We proceed with the action of the  term  $\cal{E}(\varphi_{0}):=-\sum_{j=1}^{n}e_{j}\cdot \Big[\nabla^{s}_{\nabla^{s}_{e_{j}}X}\varphi_{0}+4s\nabla^{s}_{T(X, e_{j})}\varphi_{0}\Big]$  on   $\varphi_{0}$. Because 
         \[
              \nabla^{s}_{\nabla^{s}_{e_{j}}X}\varphi_{0}=\frac{4s-1}{4}\big((\nabla^{s}_{e_{j}}X)\lrcorner T\big)\cdot\varphi_{0}, \quad  \nabla^{s}_{T(X, e_{j})}\varphi_{0}=\frac{4s-1}{4}\big(T(X, e_{j})\lrcorner T\big)\cdot \varphi_{0},
              \]
and  $\sum_{j}e_{j}\cdot \big(T(X, e_{j})\lrcorner T\big)=-\sum_{j}T(X, e_{j})\cdot (e_{j}\lrcorner T)$ (see \cite[p.~325]{ABK}), we conclude that 
 \begin{eqnarray}
          \cal{E}(\varphi_{0})&=&-\frac{(4s-1)}{4}\sum_{j}e_{j}\cdot \big((\nabla_{e_{j}}^{s}X)\lrcorner T\big)\cdot\varphi_{0} -s(4s-1)\sum_{j}e_{j}\cdot \big(T(X, e_{j})\lrcorner T\big)\cdot\varphi_{0}\nonumber \\ 
          &=&-\frac{(4s-1)}{4}\sum_{j}e_{j}\cdot \big((\nabla_{e_{j}}^{s}X)\lrcorner T\big)\cdot\varphi_{0} +s(4s-1)\sum_{j}T(X, e_{j})\cdot (e_{j}\lrcorner T)\cdot\varphi_{0}. \label{citedr} 
          \end{eqnarray}
     Now,  by the definition of $d^{s}$ and since $\nabla^{s}_{X}(Y\lrcorner T)=(\nabla^{s}_{X}Y)\lrcorner T+Y\lrcorner (\nabla^{s}_{X}T)$ for any $X, Y\in\Gamma(TM)$ and $s\in\bb{R}$, it also follows that
 \[
      d^{s}(X\lrcorner T)=\sum_{j}e_{j}\wedge \big(\nabla^{s}_{e_{j}}(X\lrcorner T)\big)=\sum_{j}e_{j}\wedge  \big((\nabla^{s}_{e_{j}}X)\lrcorner T\big)+\sum_{j}e_{j}\wedge \big(X\lrcorner (\nabla^{s}_{e_{j}}T)\big).
\]
Having in mind the isomorphism $X\cdot \simeq X\wedge-X\lrcorner$, we   combine the last relation with a part of (\ref{citedr}), i.e.    \begin{eqnarray}
\cal{A}(\varphi_{0})&:=& \frac{(4s-1)}{4}   d^{s}(X\lrcorner T)\cdot\varphi_{0}-\sum_{j=1}^{n}e_{j}\cdot \nabla^{s}_{\nabla^{s}_{e_{j}}X}\varphi_{0}\nonumber\\
&=& \frac{(4s-1)}{4}\sum_{j}\Big[e_{j}\lrcorner \big((\nabla_{e_{j}}^{s}X)\lrcorner T\big)\Big]\cdot\varphi_{0}+ \frac{(4s-1)}{4}\sum_{j}\Big[e_{j}\wedge \big(X\lrcorner (\nabla^{s}_{e_{j}}T)\big)\Big]\cdot\varphi_{0}.\label{form4}
   \end{eqnarray}
On the other hand, for any $X\in\Gamma(TM)$, $T\in\Lambda^{3}T^{*}M$  and $s\in\bb{R}$ it holds that
    \begin{eqnarray*}
             \delta^{s}(X\lrcorner T)&=&-\sum_{j}e_{j}\lrcorner \nabla^{s}_{e_{j}}(X\lrcorner T)=-\sum_{j} e_{j}\lrcorner \big((\nabla_{e_{j}}^{s}X)\lrcorner T\big)-\sum_{j}e_{j}\lrcorner \big(X\lrcorner (\nabla^{s}_{e_{j}}T)\big).
      \end{eqnarray*}  
  Therefore, adding appropriately with (\ref{form4}),    the first sums cancel   each other and we obtain 
     \begin{eqnarray}
 \cal{A}(\varphi_{0})+\frac{(4s-1)}{4}  \delta^{s}(X\lrcorner T)\cdot\varphi_{0}&:=&\frac{(4s-1)}{4}\Big[d^{s}(X\lrcorner T)+\delta^{s}(X\lrcorner T)\Big]\cdot\varphi_{0}-\sum_{j=1}^{n}e_{j}\cdot \nabla^{s}_{\nabla^{s}_{e_{j}}X}\varphi_{0}\nonumber\\
   &=& \frac{(4s-1)}{4}\sum_{j}\Big[e_{j}\wedge \big(X\lrcorner (\nabla^{s}_{e_{j}}T)\big)\Big]\cdot\varphi_{0}\nonumber\\
   &&- \frac{(4s-1)}{4}\sum_{j}\Big[e_{j}\lrcorner \big(X\lrcorner (\nabla^{s}_{e_{j}}T)\big)\Big]\cdot\varphi_{0}\nonumber\\
   &=& \frac{(4s-1)}{4}\sum_{j}e_{j}\cdot \big(X\lrcorner (\nabla_{e_{j}}^{s}T)\big)\cdot\varphi_{0},\label{form5}
   \end{eqnarray}   
   where for the last equality we apply again  the isomorphism $X\cdot \simeq X\wedge-X\lrcorner$.  In this way and by combining   the relations    (\ref{form3}), (\ref{citedr}), (\ref{form4}), (\ref{form5}) and (\ref{gen1/2ric}), we conclude that   
            \begin{eqnarray}
    \frac{1}{2}\Ric^{s}(X)\cdot\varphi_{0}&=&\frac{(4s-1)}{4}\sum_{j}e_{j}\cdot \big(X\lrcorner (\nabla_{e_{j}}^{s}T)\big)\cdot\varphi_{0} + \frac{(16s^{2}-1)}{8}\sum_{j}T(X, e_{j})\cdot (e_{j}\lrcorner T)\cdot\varphi_{0} \nonumber \\
    &&+s(3-4s)(X\lrcorner \sigma_{T})\cdot\varphi_{0}. \label{form1}
\end{eqnarray}

   \noindent The final step is based on the fact that  under the condition  $\nabla^{c}T=0$, it holds that (see \cite[Thm.~B.1]{ABK})
   \[
   (\nabla^{s}_{X}T)(Y, Z, W)=\frac{4s-1}{2}\sigma_{T}(Y, Z, W, X)=-\frac{4s-1}{2}\sigma_{T}(X, Y, Z, W).
   \]
       Thus, for any $X\in\Gamma(TM)$ the 3-form $(\nabla^{s}_{X}T)$ equals to 
   \[
   (\nabla^{s}_{X}T)=-\frac{4s-1}{2}(X\lrcorner\sigma_{T})
   \]
  and this has as a result  the following: 
   \begin{eqnarray*}
   \frac{(4s-1)}{4}\sum_{j}e_{j}\cdot \big(X\lrcorner (\nabla_{e_{j}}^{s}T)\big)&=&-\frac{(4s-1)^{2}}{8}\sum_{j}e_{j}\cdot\big(X\lrcorner (e_{j}\lrcorner\sigma_{T})\big)\cdot\varphi_{0}\\
   &=&\frac{(4s-1)^{2}}{8}\sum_{j}e_{j} \cdot\big(e_{j}\lrcorner (X\lrcorner\sigma_{T})\big)\cdot\varphi_{0}=\frac{(4s-1)^{2}}{8}\sum_{j}e_{j} \cdot \big(e_{j}\lrcorner \sigma^{X}_{T}\big)\cdot\varphi_{0}\\
   &=&\frac{3(4s-1)^{2}}{8} \sigma^{X}_{T}\cdot\varphi_{0}=\frac{3(4s-1)^{2}}{8} (X\lrcorner \sigma_{T})\cdot\varphi_{0}.
   \end{eqnarray*}
Here, $\sigma^{X}_{T}:=X\lrcorner\sigma_{T}$ is a 3-form, see also Lemma \ref{USEF3}. In combination with  (\ref{form1}), this observation  completes  the proof. $\blacksquare$ 
   
   % and finish the proof. $\blacksquare$

     \begin{remark}\label{reform}\textnormal{ \noindent For the parameter $s=1/4$, Theorem \ref{newpar}  reduces to $\Ric^{c}(X)\cdot\varphi_{0}=(X\lrcorner \sigma_{T})\cdot\varphi_{0}$, for any $X\in\Gamma(TM)$, as it should be according to \cite{FrIv}, or our Corollary \ref{genpar}. }
\end{remark}

\noindent  For completeness, let us use  (\ref{NPAR}) to  verify   the relation between the scalar curvatures $\Sca^{s}$ and $\Sca^{g}$, namely $\Sca^{s}=\Sca^{g}-24s^{2}\|T\|^{2}$ (see \cite{AF}).  Indeed, as in the proof of Theorem \ref{parallel}, we consider a local orthonormal frame $\{e_{i}\}$ and we write $\sum_{i}e_{i}\cdot \Ric^{s}(e_{i})\cdot\varphi_{0}=-\Sca^{s}\cdot \ \varphi$. On the other hand, it is $T(e_i, e_j)=\sum_{k}T^{k}_{ij}e_{k}=\sum_{k}T(e_{i}, e_{j}, e_{k})e_{k}$ and by Theorem \ref{newpar}, for a non-trivial $\nabla^{c}$-paraller spinor $\varphi_{0}$,  we see that
\begin{eqnarray*}
\sum_{i}e_{i}\cdot \Ric^{s}(e_{i})\cdot\varphi_{0}&=&-\frac{(16s^{2}-1)}{4}\sum_{i, j}e_{i}\cdot e_{j}\cdot\big(T(e_{i}, e_{j})\lrcorner T\big)\cdot \varphi_{0}+\frac{(16s^{2}+3)}{4}\sum_{i}e_{i}\cdot (e_{i}\lrcorner \sigma_{T})\cdot\varphi_{0}\\
&\overset{(\star)}{=}&-\frac{(16s^{2}-1)}{4}\sum_{i, j, k}T_{ij}^{k}e_{i}\cdot e_{j}\cdot (e_{k}\lrcorner T)\cdot\varphi_{0}+(16s^{2}+3)\sigma_{T}\cdot\varphi_{0}\\
&\overset{(\ref{ekT})}{=}&-\frac{(16s^{2}-1)}{2}\sum_{k}(e_{k}\lrcorner T)\cdot (e_{k}\lrcorner T)\cdot\varphi_{0}+(16s^{2}+3)\sigma_{T}\cdot\varphi_{0}\\
&\overset{(\ast)}{=}&-\frac{(16s^{2}-1)}{2}\Big[2\sigma_{T}-3\|T\|^{2}\Big]\cdot\varphi_{0}+(16s^{2}+3)\sigma_{T}\cdot\varphi_{0}\\
&=&4\sigma_{T}\cdot\varphi_{0}+\frac{3(16s^{2}-1)}{2}\|T\|^{2}\cdot\varphi_{0}=-\Sca^{c}\ \cdot\varphi_{0}+\frac{3(16s^{2}-1)}{2}\|T\|^{2}\cdot\varphi_{0},
\end{eqnarray*}
where for $(\ast)$ we used the fact $\sum_{j}(e_{j}\lrcorner T)\cdot (e_{j}\lrcorner T)=2\sigma_{T}-3\|T\|^{2}$ (see \cite[p.~328]{ABK}). We deduce that
\[
\Sca^{s}\cdot \ \varphi_{0}=\Sca^{c}\cdot \ \varphi_{0}-\frac{3(16s^{2}-1)}{2}\|T\|^{2}\cdot\varphi_{0}=\Sca^{g} \cdot \ \varphi_{0}-24s^{2}\|T\|^{2}\cdot\varphi_{0}.
\]
and the assertion follows since $\varphi_{0}$ does not have zeros.

\begin{remark}
\textnormal{For these computations one  could even proceed  as follows: Based on   (\ref{usef3}), in   $(\star)$ we   replace $(e_{k}\lrcorner T)\cdot\varphi_{0}$ by $\frac{4}{4s-1}\nabla^{s}_{e_{k}}\varphi_{0}$  for any $s\neq 1/4$. Then, we use  (\ref{motiv}) to obtain (a multiple of) the   operator $\slashed{D}^{s}$, appearing  in Theorem \ref{parallel}.  For the final step one needs a description of the $\slashed{D}^{s}$-eigenvalues when this operator acts on $\nabla^{c}$-parallel spinors, which we present in Section \ref{hope2}, see  Proposition \ref{parslash}.}
\end{remark}

   \begin{corol}\label{newpar2}
 Consider a triple $(M^{n}, g, T)$ as in Theorem \ref{newpar}, admitting  a non-trivial $\nabla^{c}$-parallel spinor  $\varphi_{0}\in\cal{F}^{g}(\gamma)$ $(\gamma\in\bb{R})$. Then,  the Riemannian Ricci endomorphism acts on $\varphi_{0}$ as
  \begin{eqnarray}
    \Ric^{g}(X)\cdot\varphi_{0}&=& \frac{1}{4}\sum_{j}e_{j}\cdot \big(T(X, e_{j})\lrcorner T\big)\cdot\varphi_{0} +\frac{3}{4}(X\lrcorner \sigma_{T})\cdot\varphi_{0}\nonumber \\
 &=&  \frac{1}{8}\sum_{j}e_{j}\cdot(X\lrcorner T)\cdot (e_{j}\lrcorner T)\cdot\varphi_{0}-\frac{3\gamma}{8}\cdot (X\lrcorner T)\cdot\varphi_{0}+ \frac{3}{4}(X\lrcorner \sigma_{T})\cdot\varphi_{0}.\label{malakia}
   \end{eqnarray}
   \end{corol}
     %\vskip 0.1cm
          \noindent {\bf Proof.}
       The first expression  occurs by Theorem \ref{newpar} for $s=0$.   The second one is based on the following lemma (observe that a similar reformulation as (\ref{malakia}) applies also to Theorem \ref{newpar} and hence the eigenvalue $\gamma$ can appear   in the corresponding  expression, as well).  $\blacksquare$
       
        \begin{lemma}\label{QUEST}
         For any vector field $X\in\Gamma(TM)$,   3-form $T\in\Lambda^{3}T^{*}M$ and orthonormal frame $\{e_{j}\}$,   the following holds:      
          \[
      \sum_{j=1}^{n}T(X, e_{j})\cdot (e_{j}\lrcorner T)\equiv  \sum_{j=1}^{n}(e_{j}\lrcorner X\lrcorner T)\cdot (e_{j}\lrcorner T)=-\frac{1}{2}\sum_{j=1}^{n}e_{j}\cdot (X\lrcorner T)\cdot (e_{j}\lrcorner T)+\frac{3}{2}(X\lrcorner T)\cdot T.
          \]
                    \end{lemma}
            \vskip 0.1cm
          \noindent {\bf Proof.}  By (\ref{clif1}) we see that $e_{j}\cdot \omega-(-1)^{p}\omega\cdot e_{j}=-2(e_{j}\lrcorner \omega)$ for any $p$-form $\omega\in\Lambda^{p}T^{*}M$. Since for any $X\in\Gamma(TM)$ the quantity $\omega:=X\lrcorner T$ is a 2-form, it follows that $-\frac{1}{2}\Big[e_{j}\cdot (X\lrcorner T)-(X\lrcorner T)\cdot e_{j}\Big]=e_{j}\lrcorner X\lrcorner T.$ Consequently, recalling that $\sum_{j}e_{j}\cdot (e_{j}\lrcorner T)=3T$   (cf. \cite[p.~328]{ABK} or Lemma \ref{USEF3}), one  gets the result:
                    \begin{eqnarray*}
          \sum_{j}(e_{j}\lrcorner X\lrcorner T)\cdot (e_{j}\lrcorner T)&=&-\frac{1}{2}\sum_{j}\Big[e_{j}\cdot (X\lrcorner T)-(X\lrcorner T)\cdot e_{j}\Big]\cdot (e_{j}\lrcorner T)\\
          &=&-\frac{1}{2}\sum_{j}e_{j}\cdot (X\lrcorner T)\cdot (e_{j}\lrcorner T)+\frac{1}{2}\sum_{j}(X\lrcorner T)\cdot e_{j}\cdot (e_{j}\lrcorner T)\\
          &=&-\frac{1}{2}\sum_{j}e_{j}\cdot (X\lrcorner T)\cdot (e_{j}\lrcorner T)+\frac{3}{2}(X\lrcorner T)\cdot T. \quad \blacksquare
          \end{eqnarray*}

\begin{remark}
\textnormal{As we explained in Remark \ref{ricci},  for $\nabla^{c}$-parallel torsion $T$, the Ricci tensor satisfies the relation $\Ric^{s}(X, Y)=\Ric^{g}(X, Y)-4s^{2}S(X, Y)$, where $S$ is a symmetric covariant 2-tensor defined by $S(X, Y):=\sum_{i}g\big(T(X, e_{i}), T(Y, e_{i})\big)$, see for example \cite[p.~110]{Chrysk2}. Hence,   the Ricci endomorphism $\Ric^{s}(X)$ is given by   $\Ric^{s}(X)=\Ric^{g}(X)-4s^{2}S(X)$, where $S(X)$ is the symmetric endomorphism  associated to $S$, i.e. $g(S(X), Y)=S(X, Y)$, for any $X, Y\in\Gamma(TM)$. In particular,
$\Ric^{g}(X)=\Ric^{c}(X)+\frac{1}{4}S(X)$ (cf. \cite{FrIv}). Therefore, a direct combination of Corollaries \ref{newpar2}  and  \ref{genpar} for example, allows us to describe the explicit action of $S(X)$ $(X\in\Gamma(TM))$ on $\nabla^{c}$-parallel spinors, for any metric connection $\nabla^{c}=\nabla^{g}+\frac{1}{2}T$ with $\nabla^{c}T=0$.}
\end{remark}
         
         \begin{corol}\label{Ssym}
          Consider a triple $(M^{n}, g, T)$ as in Theorem \ref{newpar}, admitting  a non-trivial $\nabla^{c}$-parallel spinor  $\varphi_{0}\in\cal{F}^{g}(\gamma)$ $(\gamma\in\bb{R})$. Then, for any $X\in\Gamma(TM)$, the action of the symmetric endomorphism $S(X)$ on $\varphi_{0}$ is given by
          \begin{eqnarray*}
          S(X)\cdot\varphi_{0}&=&\sum_{j=1}^{n}e_{j}\cdot\big(T(X, e_{j})\lrcorner T)\cdot\varphi_{0}-(X\lrcorner \sigma_{T})\cdot\varphi_{0}\\
          &=&\frac{1}{2}\sum_{j=1}^{n}e_{j}\cdot (X\lrcorner T)\cdot (e_{j}\lrcorner T)\cdot\varphi_{0}-\frac{3\gamma}{2}(X\lrcorner T)\cdot \varphi_{0}-(X\lrcorner \sigma_{T})\cdot\varphi_{0}.
          \end{eqnarray*}
         \end{corol}

\vskip 0.2cm
          \noindent   Theorem \ref{newpar} and Corollaries \ref{newpar2} and \ref{Ssym} can be applied on any triple $(M^{n}, g, T)$  endowed  with a non-integrable $G$-structure $\cal{R}\subset P^{g}$  $(G\subsetneq\SO_{n})$ and  a  $\nabla^{c}$-parallel spinor $\varphi_{0}$  with respect to the adapted (unique) characteristic connection $\nabla^{c}=\nabla^{g}+\frac{1}{2}T$, under the assumption $\nabla^{c}T=0$.    Special structures fitting in this setting are plentiful, e.g. Sasakian manifolds in any odd dimension  \cite{FrIv, FrIv2}, almost hermitian structures in even dimensions \cite{AFS, Schoe}, co-calibrated $\G_2$-structures  in dimension $7$  \cite{FrIv, FriedG2}, (non-parallel) $\Spin_7$-structures in dimension 8 \cite{Puhle, Ivspin}, to name some of them.    In the following,  we are going to illustrate   our integrability conditions   on nearly parallel $\G_2$-structures, nearly K\"ahler structures and Sasakian structures. 
          
            \subsection{Real Killing spinors which are $\nabla^{c}$-parallel}\label{pKsT}   For the first two special structures mentioned above, the description can be globalized and this is because on these manifolds %the nearly K\"ahler manifolds in dimension $6$ and the nearly parallel $\G_2$-manifolds in dimension 7, together with the 3-sphere $(\Ss^{3}, g_{\rm can})$, exhaust all compact Riemannian manifolds whose 
      the existent $\nabla^{c}$-parallel spinors coincide with the real Killling spinors, i.e. they satisfy  the additional  equation $\nabla^{g}_{X}\varphi_{0}=\kappa X\cdot\varphi_{0}$  for any $X\in\Gamma(TM)$ and some $\kappa\in\bb{R}^{*}$ (with respect to the same metric $g$ that holds $\nabla^{c}\varphi_{0}=0$).        Friedrich and Ivanov \cite[Thm.~5.6, 10.8]{FrIv} were the  first who  provided this identification and moreover  proved that any such  Einstein manifold is also $\nabla^{c}$-Einstein.    In \cite[Prop.~5.1]{Chrysk2} we generalise these results by showing  that  the Ricci endomorphism $\Ric^{s}(X)$  on such Einstein manifolds  is a multiple of the identity operator for any $s\in\bb{R}$, i.e. $  \Ric^{s}=\frac{\Sca^{s}}{n}\Id$, and moreover that the existent $\nabla^{c}$-parallel spinors are  Killing spinor with torsion with respect to $\nabla^{s}$ (or twistor spinors with torsion) for any $s\in\bb{R}\backslash\{0, 1/4\}$ (for $s=5/12$  and 6-dimensional nearly K\"ahler manifolds this result was known by \cite[Thm.~6.1]{ABK}).    A   direct and very simplified   proof of the first conclusion arises now in terms of Theorem \ref{newpar}, as follows.
   
   \vskip 0.2cm
   \noindent 
 Let us   denote by $\cal{K}^{s}(M, g)_{\zeta}:=\{\varphi\in\cal{F}^{g}:  \nabla^{s}_{X}\varphi=\zeta X\cdot\varphi \ \ \forall X\in\Gamma(TM)\}$   the set of all  {\it  Killing spinors with torsion}  (KsT in short), with respect to the family $\nabla^{s}=\nabla^{g}+2sT$ $(s\neq 0)$,   with   Killing number    $\zeta\neq 0$ (we refer to \cite{ABK, Chrysk2} for a detailed exposition related to this kind of spinors). Similarly, we shall write  $\cal{K}^{0}(M^{n}, g)_{\kappa}$ for the set of all real Killing spinors with Killing number $\kappa\neq 0$.   Assume   that   $(M^{n}, g, T)$  is a compact connected Riemannian spin manifold $(M^{n}, g, T)$, with $\nabla^{c}T=0$ and   positive scalar curvature given by $\Sca^{g}=\frac{9(n-1)\gamma^{2}}{4n}$, for some constant $0\neq\gamma\in\Spec(T)$.   In  \cite[Thm.~3.7]{Chrysk2} we extended the    identification mentioned above, namely 
     \[
  \ke\big(\nabla^{c}\big) \cong  \bigoplus_{\gamma\in\Spec(T)}\Big[\Gamma(\Sigma_{\gamma})\cap \cal{K}^{0}(M^{n}, g)_{\frac{3\gamma}{4n}}\Big],
  \]
  by proving that   for any non-trivial  $\nabla^{c}$-parallel spinor $\varphi_{0}$    the following  conditions are equivalent:
    \begin{itemize}
\item[$(a)$] $\varphi_{0}\in\cal{F}^{g}({\gamma})\cap\ke(\cal{P}^{s}):=\ke(\cal{P}^{s}\big|_{\Sigma_{\gamma}^{g}M})$  with respect to the family $\{\nabla^{s} : s\in\bb{R}\backslash\{1/4\}\}$, \\
\item[$(b)$] $\varphi_{0}\in \cal{K}^{s}(M, g)_{\zeta}$  with respect to  the family $\{\nabla^{s} : s\in\bb{R}\backslash\{0, 1/4\}\}$ with $\zeta:=3(1-4s)\gamma/4n$,\\
\item[$(c)$]  $ \varphi_{0}\in \cal{K}^{0}(M, g)_{\kappa}$ with $\kappa:={3\gamma}/{4n}$.
 \end{itemize}
 This correspondence  allows us   now to proceed with the following (see \cite[Prop.~5.1]{Chrysk2} for another method): %The special cases  $s=0$ and $s=1/4$ were known before, by  \cite{Gr, FKMS}  and   \cite{FrIv}, respectively, while notice that for $n=3$, the manifold must be $\Ric^{c}$-flat and hence isometric to the  3-sphere $(\Ss^{3}, g_{\rm can})$.  Now we are ready to prove% For instance, this applies  to  6-dimensional nearly K\"ahler manifolds, or  to 7-dimensional nearly-parallel $\G_2$-manifolds, where the $\nabla^{c}$-parallel spinors coincide with the real Killing spinors, see also Examples \ref{npg2}, \ref{nKm} and for more details we refer to \cite{FrIv, Chrysk2}. % and \cite{FrIv, Chrysk2}. 

  \begin{corol}\label{mine22}
 Let $(M^{n}, g, T)$ $(n\geq 3)$ be a compact Riemmanian spin manifold,  endowed with the characteristic connection   $\nabla^{c}=\nabla^{g}+\frac{1}{2}T$, such that $\nabla^{c}T=0$.  Assume that $0\neq \varphi_{0}\in\Sigma^{g}_{\gamma}M$ $(\bb{R}\ni\gamma\neq 0)$   is a non-trivial $\nabla^{c}$-parallel spinor, which satisfies one of the conditions $(a)$,  $(b)$, or $(c)$. Then, the $\frac{1}{2}$-$\Ric^{s}$-formula gives rise to the equation
 \[
      \Ric^{s}(X)\cdot\varphi_{0}=\frac{\Sca^{s}}{n}X\cdot\varphi_{0}=\frac{3\gamma^{2}(-3 + 3n - 144s^{2} + 16ns^{2})}{4n^{2}}X\cdot\varphi_{0}, \quad     (\eighthnote)
\]
for any $X\in\Gamma(TM)$, where $\Sca^{s}:=\frac{6\gamma^{2}}{n}\Big[\frac{6(n-1)(1-4s)^{2}+96s(1-4s)+16s(3-4s)(n-3)}{16}\Big]=\frac{3\gamma^{2}(-3 + 3n - 144s^{2} + 16ns^{2})}{4n}$. Moreover, for $n\neq 9$   the symmetric endomorphism $S(X)$ acts on $\varphi_{0}$ as a multiple of the identity,  
\begin{equation}\label{actionS}
S(X)\cdot\varphi_{0}=-\frac{3\gamma^{2}(n-9)}{n^{2}}X\cdot\varphi_{0}.
\end{equation}
  \end{corol}
  
      \vskip 0.1cm
          \noindent {\bf Proof}. 
           Assume that $0\neq \varphi_{0}\in\Sigma^{g}_{\gamma}M$  is a $\nabla^{c}$-parallel spinor which satisfies any of the conditions  $(a)$, $(b)$, or $(c)$.   Then, by  \cite[Prop.~3.2, Thm.~3.7]{Chrysk2}    this is equivalent to say that $\varphi_{0}$ is a solution of the equation
           \begin{equation}\label{t1s}
(X\lrcorner T)\cdot \varphi_{0}+\frac{3\gamma}{n}X\cdot\varphi_{0}=0,
\end{equation} 
for any $X\in\Gamma(TM)$.  Moreover, the Ricci tensor $\Ric^{c}$ is computed algebraically (see for example the proof of \cite[Prop.~5.1, (b)]{Chrysk2})
\begin{equation}\label{ricccein}
\Ric^{c}(X)\cdot\varphi_{0}=(X\lrcorner \sigma_{T})\cdot\varphi_{0}=\frac{3\gamma^{2}(n-3)}{n^{2}}X\cdot\varphi_{0}.
\end{equation}
By (\ref{t1s}) it follows that an arbitrary vector field $X$  satisfies the equation   $(T(X, e_{j})\lrcorner T)\cdot\varphi_{0}=-\frac{3\gamma}{n}T(X, e_{j})\cdot\varphi_{0}$ and because $\sum_{j} e_{j}\cdot T(X, e_{j})=2(X\lrcorner  T)$,  an application  of Theorem \ref{newpar} gives rise to
\[
\Ric^{s}(X)\cdot\varphi_{0}=-\frac{18\gamma^{2}(16s^{2}-1)}{4n^{2}}X \cdot\varphi_{0}+\frac{3\gamma^{2}(n-3)(16s^{2}+3)}{4n^{2}}X\cdot\varphi_{0},
\]
  which equals  to the given relation. Finally, Corollary \ref{Ssym} yields the  expression for the action of $S(X)$.   $\blacksquare$ 

%\begin{eqnarray*}
%&=&\frac{3\gamma(16s^{2}-1)}{4n}\sum_{j}e_{j}\cdot T(X, e_{j})\cdot\varphi_{0}+\frac{3\gamma^{2}(n-3)(16s^{2}+3)}{4n^{2}}X\cdot\varphi_{0}\\
%&=&\frac{6\gamma(16s^{2}-1)}{4n} (X\lrcorner T)\cdot\varphi_{0}+\frac{3\gamma^{2}(n-3)(16s^{2}+3)}{4n^{2}}X\cdot\varphi_{0}\\
%\end{eqnarray*}

   \vskip 0.2cm
     \noindent  Therefore, as in \cite{Chrysk2}, one concludes that a (complete) Riemannian manifold $(M^{n}, g, T)$ $(n>3)$ satisfying Corollary \ref{mine22}   is a compact Einstein  manifold with constant positive scalar curvature $\Sca^{g}=\frac{9(n-1)\gamma^{2}}{4n}$,  a $\nabla^{c}$-Einstein manifold with parallel torsion and constant positive scalar curvature $\Sca^{c}=\frac{3(n-3)\gamma^{2}}{n}$ and  moreover that satisfies  the ``harmony equation''  $(\eighthnote)$ for any other $s$.  For  $n=3$,    $(M^{3}, g, T)$  is $\Ric^{c}$-flat and hence isometric to the 3-sphere. Notice  also  that under our point of view,  the study of 6-dimensional nearly K\"ahler manifolds and  7-dimensional nearly parallel $\G_2$-manifolds reduces to be qualitatively the same. Therefore, below  we illustrate our conclusions  only for one of these two classes, e.g. nearly parallel $\G_2$-manifolds and  similarly is treated the former class  (see also Example \ref{nKm} and \cite{FrIv, ABK, Chrysk2} for useful details).
     
       %\subsection{Nearly parallel $\G_2$-structures}\label{npg2}
        %Consider a 7-dimensional oriented Riemannian manifold $(M^7,g)$. We say that $M$ admits  a $\G_2$-structure if the $\SO_{7}$-principal bundle $P^{g}$ of positively oriented orthonormal frames admits a reduction to $\G_2\subset\SO_7$. Alternatively, this amounts to the choice of a generic 3-form $\omega$, which in terms of     a local orthonormal frame $\{e_{1}, \ldots, e_{7}\}$ may  be expressed by
              %    \[
      %   \omega:=e_{127}+e_{135} -e_{146}-e_{236}-e_{245} +e_{347} +e_{567},
      %   \]
       %  where in general by $e_{i_{1}\ldots i_{p}}$ we denote  the $p$-form $e_{i_{1}}\wedge \ldots \wedge e_{i_{p}}$.  A $\G_2$-structure, say $P_{\G_2}\subset P^{g}$, induces a spin structure $\widetilde{P}^{g}:=P_{\G_2}\times_{\G_2}\Spin_{7}$ and defines a nonwhere vanishing spinor $\varphi_{0}\in\Gamma(\Sigma^{g}M):=\cal{F}^{g}$, where  $\Sigma^{g}M=P_{\G_2}\times_{\G_2}\Delta_{7}$ is the spinor bundle. In fact, because $\G_{2}$ preserves both $\omega$ and   $\varphi_{0}\in\Delta_{7}$,  they   induce necessarily    the same data, namely   
  %$    \omega(X, Y, Z):=(X\cdot Y\cdot Z\cdot\varphi_{0}, \varphi_{0})$.  
  \begin{example}
   \textnormal{Consider a  nearly parallel $\G_2$-manifold $(M^{7}, g, \omega)$, i.e. a  7-dimensional oriented Riemannian manifold with   a $\G_2$-structure $\omega\in\Gamma(\Lambda^{3}_{+}T^{*}M)$ satisfying the differential equation $  d\omega=-\tau_{0}\ast \omega$,   for some real constant $\tau_{0}\neq 0$ (we refer to \cite{FKMS, FrIv, FriedG2} for an introduction to $\G_2$-structures and also to $\G_2$-structures carrying  a characteristic connection). In \cite[Cor.~4.9]{FrIv}  it was shown  that a nearly parallel $\G_2$-manifold admits a unique characteristic connection $\nabla^{c}$ with {\it parallel} skew-torsion $T$, given by $T:=\frac{1}{6}(d\omega, \ast\omega)\cdot\omega$, in particular $T=-\frac{\tau_{0}}{6}\omega$ and  $\|T\|^{2}=\frac{7}{36}\tau_{0}^{2}$. Moreover,  there exists a unique spinor field  $\varphi_{0}$ which  is  $\nabla^{c}$-parallel (cf. \cite[Prop.~3.2]{FriedG2})  and satisfies the equation $T\cdot \varphi_{0}=-\frac{7\tau_{0}}{6}\varphi_{0}=-\sqrt{7}\|T\|\varphi_{0}$, i.e. $\gamma=-\sqrt{7}\|T\|$.\footnote{Notice that here our 3-form $\omega$ is such that $\omega\cdot\varphi_{0}=7\varphi_{0}$, see  \cite[Lem.~2.3]{ACFH}.}  In fact, $\varphi_{0}$ is a  real Killing spinor and hence $\Ric^{g}(X)\cdot\varphi_{0}=\frac{27}{14}\|T\|^{2}X\cdot\varphi_{0}=\frac{3\tau_{0}^{2}}{8}X\cdot\varphi_{0}$  (\cite{FKMS,   FrIv}) . More general,  in  \cite[Exam.~5.3]{Chrysk2}  we  deduced the ``harmony equation'' $(\eighthnote)$, i.e.  
   \[       \Ric^{s}(X)\cdot\varphi_{0}=\frac{3(9-16s^{2})}{14}\|T\|^{2}X\cdot\varphi_{0}, % \quad (\natural)
       \]
 for any $s\in\bb{R}$ and $X\in\Gamma(TM)$.  Let us  provide a new proof of this result,   via Theorem  \ref{newpar}.  As in the proof of Corollary \ref{mine22}, the key point  is that  $\varphi_{0}$ satisfies  the  equations (\ref{t1s}) and (\ref{ricccein}) respectively, i.e. 
          \[
 (X\lrcorner T)\cdot\varphi_{0}=\frac{\tau_{0}}{2}X\cdot\varphi_{0}=\frac{3\|T\|}{\sqrt{7}}X\cdot\varphi_{0}, \quad \Ric^{c}(X)\cdot\varphi_{0}=(X\lrcorner \sigma_{T})\cdot\varphi_{0}=\frac{12}{7}\|T\|^{2}X\cdot\varphi_{0},
 \]
 for any $X\in\Gamma(TM)$.   Both of them can be found  \cite{Chrysk2} (see also \cite[Lem.~2.3]{ACFH} and  \cite[p.~318]{FrIv}).  The first represents the Killing equation, or equivalent the twistor equation \cite[Prop.~3.2, Thm.~4.2]{Chrysk2}, while   the second one states that $(M^{7}, g, \omega)$ is a $\nabla^{c}$-Einstein manifold. Hence, (\ref{NPAR}) yields that
 \begin{eqnarray*}
 \Ric^{s}(X)\cdot\varphi_{0}&=&-\frac{3(16s^{2}-1)\|T\|}{4\sqrt{7}}\sum_{j=1}^{7}e_{j}\cdot T(X, e_{j})\cdot\varphi_{0}+\frac{12(16s^{2}+3)\|T\|^{2}}{28}X\cdot\varphi_{0}\\
 &=&-\frac{6(16s^{2}-1)\|T\|}{4\sqrt{7}} (X\lrcorner T)\cdot\varphi_{0}+\frac{12(16s^{2}+3)\|T\|^{2}}{28}X\cdot\varphi_{0}\\
 &=&-\frac{18(16s^{2}-1)\|T\|^{2}}{28} X\cdot\varphi_{0} +\frac{12(16s^{2}+3)\|T\|^{2}}{28}X\cdot\varphi_{0}, \end{eqnarray*}
 which gives rises to the result. Finally, by (\ref{actionS}) we compute $S(X)\cdot\varphi_{0}=\frac{6\|T\|^{2}}{7}X\cdot\varphi_{0}$ (cf. \cite{FrIv}).}
\end{example}

         \subsection{5-dimensional Sasakian structures}\label{sasak}
          \textnormal{Recall that a Sasakian structure on a Riemannian manifold $(M^{2n+1}, g)$ consists of a Killing vector field $\xi$ of unit length, the so-called Reeb vector field, such that the endomorphism $\phi : TM\to TM$ given by $\phi(X)=-\nabla^{g}_{X}\xi$, satisfies $(\nabla^{g}_{X}\phi)(Y)=g(X, Y)\xi-g(\xi, Y)X$ for any $X, Y\in\Gamma(TM)$. The dual 1-form $\eta$ of $\xi$  solves the equation $d\eta=2F$,  where $F(X, Y):=g(X, \phi(Y))$ is the fundamental 2-form, see  for example \cite{FKim,  Boy1} for equivalent definitions and more details.   Let us focus on   5-dimensional Sasakian manifolds $(M^{5}, g, \xi, \eta, \phi)$. We fix an orthonormal basis $e_1, \ldots, e_5$ of $T_{x}M\simeq\bb{R}^{5}$ and use the abbreviation $e_{i_{1}\ldots i_{p}}$  for the $p$-form $e_{i_{1}}\wedge\ldots\wedge e_{i_{p}}$.  It is
          \[
          \xi:=e_5, \quad \phi:=-(e_{12}+e_{34}), \quad F:=e_{12}+e_{34},
          \]
           and in terms of $\phi$, our orthonormal frame reads by  $\{e_{1}, e_{2}:=-\phi(e_{1}), e_{3}, e_{4}=-\phi(e_{3}), e_{5}=\xi\}$, with  $\phi(\xi)=0$.    By \cite[Prop.~7.1]{FrIv} it is known that there exists a unique metric connection $\nabla^{c}$ with {\it parallel} skew-torsion
          \[
          T=\eta\wedge d\eta=2\eta\wedge F=2(e_{125}+e_{345}),
          \]
                     preserving the Sasakian structure, $\nabla^{c}g=\nabla^{c}\eta=\nabla^{c}\phi=0$.  The torsion form  $T$  acts on the 5-dimensional spin representation $\Delta_{5}$ 
        with eigenvalues $(-4, 0, 0, 4)$. Hence, the spinor bundle $\Sigma^{g}M$ splits into two 1-dimensional subbundles  and one 2-dimensional subbundle, i.e. $  \Sigma^{g}M=\Sigma^{g}_{-4}M\oplus \Sigma^{g}_{0}M\oplus\Sigma_{4}^{g}M$ with $\Sigma_{\pm 4}^{g}M:=\{\varphi\in\Sigma^{g}M : T\cdot\varphi=\pm 4\varphi\}$ and $\Sigma^{g}_{0}M:=\{\varphi\in\Sigma^{g}M : T\cdot\varphi=0\}$, respectively.}
        
        \vskip 0.2cm
              \noindent\textnormal{In the direction of the Reeb vector field $\xi$  the Riemannian Ricci endomorphism must occur with eigenvalue $4$, i.e. $\Ric^{g}(\xi)=4\xi$ (cf. \cite{FKim}).  Let us explain how  Corollary \ref{newpar2} fits with this result.   Assume  that there exists some $\nabla^{c}$-parallel spinor $\varphi_{1}$, which for instance belongs to     $\Sigma^{g}_{-4}M$, i.e. $T\cdot\varphi_{1}  =-4\varphi_{1}$.   Any vector field $X$ satisfies  $X\lrcorner F=-\phi(X)$, hence  by (\ref{clif1}) we get that (see also \cite[Lem.~6.3]{FKim})
        \begin{equation}\label{contact}
        X\cdot d\eta-d\eta\cdot X=-2(X\lrcorner d\eta)=-4(X\lrcorner F)=4\phi(X), \quad \forall X\in\Gamma(TM).
        \end{equation}
It is $\sigma_{T}=4e_{1234}$, $\xi\lrcorner \sigma_{T}=0$   and $\xi\lrcorner T=d\eta$.  Thus,  by applying for example (\ref{malakia}) (for a description based on the first expression of Corollary \ref{newpar2}, see the proof of Theorem \ref{genparsak} below), we obtain
        \begin{eqnarray*}
        \Ric^{g}(\xi)\cdot\varphi_{1}%&=&%\frac{1}{8}\sum_{j}e_{j}\cdot(\xi\lrcorner T)\cdot (e_{j}\lrcorner T)\cdot\varphi_{1}-\frac{3}{8}\gamma\cdot (\xi\lrcorner T)\cdot\varphi_{1}+ \frac{3}{4}(\xi\lrcorner \sigma_{T})\cdot\varphi_{1}\\
&=&\frac{1}{8}\sum_{j}e_{j}\cdot d\eta\cdot (e_{j}\lrcorner T)\cdot\varphi_{1}+\frac{3}{2}d\eta\cdot\varphi_{1}\\
        &\overset{(\ref{contact})}{=}&\frac{1}{8}\sum_{j}d\eta\cdot e_{j}\cdot (e_{j}\lrcorner T)\cdot\varphi_{1}+\frac{1}{2}\sum_{j}\phi(e_{j})\cdot (e_{j}\lrcorner T)\cdot\varphi_{1}+\frac{3}{2}d\eta\cdot\varphi_{1}\\
        &=&\frac{3}{8}d\eta\cdot T\cdot\varphi_{1}+\frac{1}{2}\sum_{j}\phi(e_{j})\cdot (e_{j}\lrcorner T)\cdot\varphi_{1}+\frac{3}{2}d\eta\cdot\varphi_{1}=\frac{1}{2}\sum_{j}\phi(e_{j})\cdot (e_{j}\lrcorner T)\cdot\varphi_{1}.
        \end{eqnarray*}
       One also computes $e_{2}\lrcorner T=-2e_{15}$, $e_{1}\lrcorner T=2e_{25}$,  $e_{4}\lrcorner T= -2e_{35}$, and $e_{3}\lrcorner T=2e_{45}$, which finally yields the desired assertion: $\Ric^{g}(\xi)\cdot\varphi_{1}=\big(-e_{2}\cdot e_{25}-e_{1}\cdot e_{15}-e_{4}\cdot e_{45}-e_{3}\cdot e_{35}\big)\cdot\varphi_1=4e_{5}\cdot\varphi_1$.     %   \begin{eqnarray*}
   %  \Ric^{g}(\xi)\cdot\varphi_{1}%&=&\frac{1}{2}\Big[\phi(e_1)\cdot (e_{1}\lrcorner T)+\phi(e_{2})\cdot (e_{2}\lrcorner T)+\phi(e_{3})\cdot (e_{3}\lrcorner T)+\phi(e_{4})\cdot (e_{4}\lrcorner T)\Big]\cdot\varphi_{1}\\
   %  &=&\Big[-e_{2}\cdot e_{25}-e_{1}\cdot e_{15}-e_{4}\cdot e_{45}-e_{3}\cdot e_{35}\Big]\cdot\varphi_1\\
   %  &=&-\Big[e_{2}\cdot e_{2}\cdot e_{5}+e_{1}\cdot e_{1}\cdot e_{5}+e_{4}\cdot e_{4}\cdot e_{5}+e_{3}\cdot e_{3}\cdot e_{5}\Big]\cdot\varphi_1=4e_{5}\cdot\varphi_1.
    %   \end{eqnarray*}
   In fact, the relation $\Ric^{g}(\xi)=4\xi$ can be also obtained by using a $\nabla^{c}$-parallel spinor in  $\Sigma^{g}_{4}M$, or  in $\Sigma^{g}_{0}M$.  
 More general, in the simply-connected case  one can  use Corollary \ref{newpar2}  to verify \cite[Thm.~7.3, 7.6]{FrIv}, i.e. the fact that the existence of a $\nabla^{c}$-parallel spinor in $\Sigma^{g}_{\pm 4}M$ (resp. $\Sigma^{g}_{0}M$), requires  that $(6, 6, 6, 6, 4)$ (resp. $(-2, -2, -2, -2, 4)$) are the eigenvalues  of the Ricci endomorphism $\Ric^{g}(X)$, and the converse.   
     Next we are going to extend these known results, for any  $s\in\bb{R}$, via Theorem \ref{newpar}.  Notice that simply-connected  Sasakian spin manifolds  $(M^{5}, g, \xi, \eta, \phi)$ whose Ricci tensors $\Ric^{g}$ and $\Ric^{c}$ satisfy   the   prescribed curvature conditions, are for instance circle bundles over 4-dimensional K\"ahler-Einstein manifolds  with positive scalar curvature and the 5-dimensional Heisenberg group (in fact, Sasakian structures  $(M^{5}, g, \xi, \eta, \phi)$  with a $\nabla^{c}$-parallel spinor in $\Sigma^{g}_{0}M$ are locally equivalent to the 5-dimensional  Heisenberg group, see \cite{FrIv, FrIv2}). These examples were described in  \cite[Examp.~7.4, 7.7]{FrIv} (see also \cite[Examp.~6.1]{FKim} and \cite{FrIv2}) and    for  these Sasakian manifolds also  the following more general theorem makes sense. }

             \begin{theorem}\label{genparsak}
    Consider a 5-dimensional simply-connected Sasakian spin manifold $(M^{5}, g, \xi, \eta, \phi)$ with its characteristic connection $\nabla^{c}=\nabla^{g}+\frac{1}{2}\eta\wedge d\eta=\nabla^{g}+\eta\wedge F$.  Then,\\ %Then, the following hold:\\
     (1) There exists a $\nabla^{c}$-parallel spinor $\varphi_{1}\in\Sigma^{g}_{-4}M$, or $\varphi_{1}\in\Sigma^{g}_{4}M$, if and only if  for any  $s\in\bb{R}$ the  eigenvalues of the  Ricci tensor $\Ric^{s}$    are given by 
     \[
     \Big\{(6-32s^{2}), (6-32s^{2}), (6-32s^{2}), (6-32s^{2}), -4(16s^{2}-1)\Big\}.
     \]
         (2) There exists a $\nabla^{c}$-parallel spinor $\varphi_{0}\in\Sigma^{g}_{0}M$,  if and only if  for any  $s\in\bb{R}$  the eigenvalues of the  Ricci tensor $\Ric^{s}$   are given by 
     \[
     \Big\{-(2+32s^{2}), -(2+32s^{2}), -(2+32s^{2}), -(2+32s^{2}), -4(16s^{2}-1)\Big\}.
     \]
          \end{theorem}  
          
          \vskip 0.2cm
\noindent {\bf Proof.}  We begin again with the action of the endomorphism $\Ric^{s}(\xi)$. As before,  this is independent of which  subbundle $\Sigma^{g}_{\gamma}M$ $(\gamma\in\{-4, 0, 4\})$ the $\nabla^{c}$-parallel spinor is lying in.  So,   assume   that $\psi$ is a $\nabla^{c}$-parallel spinor such that $\psi\in \Sigma^{g}_{\gamma}M$ for some (constant)  $\gamma\in\bb{R}$.  For the computation of the first term in Theorem \ref{newpar}, it is useful to remind that   the Reeb vector field $\xi$ is $\nabla^{c}$-parallel, in particular $\xi$ is  a Killing vector field and consequently  $\nabla^{g}_{X}\xi=\frac{1}{2}T(\xi, X)=-\phi(X)=\frac{1}{2}X\lrcorner d\eta$, see \cite{FrIv2, AF}. Thus we compute
 \[
  \nabla^{g}_{e_{1}}\xi=e_{2}, \quad \nabla^{g}_{e_{2}}\xi=-e_{1}, \quad \nabla^{g}_{e_{3}}\xi=e_{4}, \quad \nabla^{g}_{e_{4}}\xi=-e_{3}, \quad \nabla^{g}_{\xi}\xi=0.
  \]
  Let us also set $\cal{W}:=-\sum_{j=1}^{5}e_{j}\cdot \big(T(\xi, e_{j})\lrcorner T\big)=\sum_{j=1}^{5}T(\xi, e_{j})\cdot (e_{j}\lrcorner T)$. Then we deduce that
 \begin{eqnarray*}
 \cal{W}&=&2\Big[(\nabla^{g}_{e_{1}}\xi)\cdot (e_1\lrcorner T)+(\nabla^{g}_{e_{2}}\xi)\cdot (e_2\lrcorner T)+(\nabla^{g}_{e_{3}}\xi)\cdot (e_3\lrcorner T)+(\nabla^{g}_{e_{4}}\xi)\cdot (e_4\lrcorner T)\Big]\\
 &=&4\Big[e_{2}\cdot e_{25}+e_{1}\cdot e_{15}+e_{4}\cdot e_{45}+e_{3}\cdot e_{35}\Big]=-16\xi.
 \end{eqnarray*}
      One finishes with the first term, after  a multiplication with  the coefficient  $(16s^{2}-1)/4$, 
   \[
     - \frac{16s^{2}-1}{4}\sum_{j=1}^{5}e_{j}\cdot \big(T(\xi, e_{j})\lrcorner T\big)\cdot\psi=\frac{16s^{2}-1}{4}\sum_{j=1}^{5}T(X, e_{j})\cdot (e_{j}\lrcorner T)\cdot\psi=-4(16 s^{2}-1)\xi\cdot\psi.
   \]
Since $\xi\lrcorner\sigma_{T}=0$, our claim follows,  $\Ric^{s}(\xi)\cdot\psi=-4(16s^{2}-1)\xi\cdot\psi$, for any $s\in\bb{R}$.  Let us proceed now with  the action of $\Ric^{s}(X)$, for some $X\in\{e_1, \ldots, e_4\}$.  We analyse  only the case $X=e_1$  and similarly are treated the other vectors.   At this point it is sufficient to assume that $\psi\in\Sigma^{g}_{\gamma}M$ $(\gamma\in\bb{R})$ is a $\nabla^{c}$-parallel spinor (we   use the fact that  $\psi:=\varphi_1\in \Sigma^{g}_{\pm 4}M$, or $\psi:=\varphi_{0}\in\Sigma^{g}_{0}M$, only at the final step).  We compute  $T(e_{1}, e_{1})=T(e_{1}, e_{3})=T(e_{1}, e_{4})=0$ and $T(e_{1}, e_{2})=2e_{5}$, $T(e_{1}, \xi)=-T(\xi, e_1)=-2\nabla^{g}_{e_{1}}\xi=-2e_{2}$.  Hence  for the first     term in   Theorem \ref{newpar}  we deduce that
 \begin{eqnarray*}
 -\sum_{j}e_{j}\cdot \big(T(e_{1}, e_{j})\lrcorner T\big)&=&\sum_{j}T(e_{1}, e_{j})\cdot (e_{j}\lrcorner T)=\Big[T(e_{1}, e_{2})\cdot (e_{2}\lrcorner T)+T(e_{1}, \xi)\cdot (\xi\lrcorner T)\Big]\\
 &=&-2\Big[2e_{5}\cdot e_{15}+e_{2}\cdot d\eta\Big]=-4\Big[2e_1+H\Big],
  \end{eqnarray*}
 since     inside $\Cl_{5}$ we get $e_{2}\cdot d\eta=2(e_{1}+H)$, where  $H:=e_{234}=\frac{1}{4}(e_{1}\lrcorner \sigma_{T})$.      Multiplying with the coefficient $(16s^{2}-1)/4$, this gives rise to
\[
        -\frac{16s^{2}-1}{4}\sum_{j=1}^{5}e_{j}\cdot \big(T(e_{1}, e_{j})\lrcorner T\big)\cdot\psi=\frac{16s^{2}-1}{4}\sum_{j=1}^{5}T(e_{1}, e_{j})\cdot (e_{j}\lrcorner T)\cdot\psi=-(16s^{2}-1)\Big[2e_1+H\Big]\cdot\psi.
\]
       Moreover, it is  $     \frac{(16s^2 + 3)}{4}(e_1\lrcorner \sigma_{T})\cdot\psi=(16s^2 + 3)H\cdot\psi$ and thus 
  \[
    \Ric^{s}(e_1)\cdot\psi =2(1-16s^{2})e_{1}\cdot\psi+4H\cdot\psi.
    \]
   The final step  includes   the action of the  3-form $H:=e_{234}=\frac{1}{4}(e_1\lrcorner \sigma_{T})$ on $\Sigma^{g}_{\pm 4}M$ and $\Sigma^{g}_{0}M$, respectively, which of course is  related with the endomorphism $\Ric^{c}(e_1)=(e_1\lrcorner \sigma_{T})$ and is computed algebraically, see \cite[pp.~324-325]{FrIv}: 
   \[
 H\cdot\psi=\left\{
  \begin{tabular}{rcl}
$e_{1}\cdot\psi$, & if &  $\psi:=\varphi_1\in\Sigma^{g}_{\pm} 4M,$\\
$-e_{1}\cdot\psi$, & if & $\psi:=\varphi_{0}\in\Sigma^{g}_{0} 4M.$
\end{tabular}\right.
\]
Consequently 
\begin{equation}\label{newsaksak}
\Ric^{s}(e_1)\cdot\psi=
\left\{
\begin{tabular}{lcl}
$2(1-16s^{2})e_{1}\cdot\psi+4e_1\cdot\psi=(6-32s^{2})\psi,$ & if & $\psi:=\varphi_{1}\in\Sigma^{g}_{\pm 4}M$,\\
$2(1-16s^{2})e_{1}\cdot\psi-4e_1\cdot\psi=-(2+32s^{2})\psi,$ & if & $\psi:=\varphi_{0}\in\Sigma^{g}_{0}M$,
\end{tabular}\right.
\end{equation}
 for any $s\in\bb{R}$.   For the converse, assume  that  $(M^{5}, g, \xi, \eta, \phi)$ is a 5-dimensional simply-connected Sasakian spin manifold whose Ricci tensor $\Ric^{s}$ $(s\in\bb{R})$ satisfies (\ref{newsaksak}).   Then, for   $s=0, 1/4$ we see that  (\ref{newsaksak}) induces the desired  prescribed conditions, i.e. $\Ric^{g}={\rm diag}(6, 6, 6, 6, 4)$,  $\Ric^{c}={\rm diag}(4, 4, 4, 4, 0)$ for  $\varphi_{1}\in\Sigma^{g}_{\pm 4}M$ and $\Ric^{g}={\rm diag}(-2, -2, -2, -2, 4)$,  $\Ric^{c}={\rm diag}(-4, -4, -4, -4, 0)$ for $\varphi_{0}\in\Sigma^{g}_{0}M$, respectively.  Hence the assertion follows as in \cite{FrIv}. This finishes the proof. $\blacksquare$       
        
  \begin{remark}
  \textnormal{The stated expressions of the Ricci endomorphism $\Ric^{s}(X)$ can be also obtained by applying  the general type $\Ric^{s}=\Ric^{g}-4s^{2}S=\Ric^{c}-\frac{(16s^{2}-1)}{4}S$, where for the action of the symmetric endomorphism $S(X)$ on the related $\nabla^{c}$-parallel spinor $\psi\in\Sigma^{g}_{\gamma}M$   one can apply Lemma \ref{Ssym}.  We refer  also  to \cite{FrIv} for   $S(X)$.}
  \end{remark}

\section{On the  differential operator $\slashed{D}^{s}=\sum_{i}(e_{i}\lrcorner T)\cdot\nabla^{s}_{e_{i}}$}\label{hope2}
 \subsection{Special $\slashed{D}^{s}$-eigenspinors} Next  we examine some special eigenspinors of  the differential operator 
\[
\slashed{D}^{s}(\varphi)=\sum_{i}(e_{i}\lrcorner T)\cdot\nabla^{s}_{e_{i}}\varphi=\slashed{D}^{0}(\varphi)+s
\sum_{i}(e_{i}\lrcorner T)\cdot(e_{i}\lrcorner T)\cdot\varphi=\slashed{D}^{0}(\varphi)+s\cal{T}\cdot\varphi_{0},
\]
appearing  in Theorem \ref{parallel},   see \cite{Agr03, AF}. Here, $\slashed{D}^{0}$   denotes the part corresponding to  the Riemannian connection $\nabla^{0}\equiv\nabla^{g}$ and  $\cal{T}:=\sum_{j} (e_{j}\lrcorner T)\cdot (e_{j}\lrcorner T)=2\sigma_{T}-3\|T\|^{2}$.  Notice also by   (\ref{motiv}), that
\[
\slashed{D}^{s}(\varphi)=\frac{1}{2}\sum_{i, j}e_{i}\cdot e_{j}\cdot\nabla^{s}_{T(e_{i}, e_{j})}\varphi.
\]
In fact, this  formula   holds also in the general case (although in the proof of Theorem \ref{parallel} we use the assumption $\nabla^{c}T=0$, this does not effect to the computations related to $\slashed{D}^{s}$).   By Proposition \ref{usef1}, (5) and Remark \ref{ricci}, one  also has  (see \cite{FrIv, AF, ABK})
\begin{equation}\label{gens1}
 \slashed{D}^{s}(\varphi)=-\frac{1}{2}\Big[D^{s}(T\cdot\varphi)+T\cdot D^{s}(\varphi)-(dT+\delta T)\cdot\varphi+8s\sigma_{T}\cdot\varphi\Big].
\end{equation}
  Let us focus  now on  triples $(M^{n}, g, T)$ with $\nabla^{c}$-parallel skew-torsion,  $\nabla^{c}T=0$.  In this case the operator $\slashed{D}^{s}$ has more equivalent expressions.
\begin{lemma}\label{equivex} \textnormal{(\cite{AF})}
Consider   a Riemannian spin manifold $(M^{n}, g, T)$ $(n\geq 3)$  endowed with a non-trivial 3-form $T\in\Lambda^{3}T^{*}M$, such that   $\nabla^{c}T=0$, where $\nabla^{c}:=\nabla^{g}+\frac{1}{2}T$. Then,  the operator $\slashed{D}^{s}$ is given by
\begin{eqnarray}
\slashed{D}^{s}(\varphi)&=&-\frac{1}{2}\Big[D^{s}(T\cdot\varphi)+T\cdot D^{s}(\varphi)-2(1-4s)\sigma_{T}\cdot\varphi\Big]\label{clasdif1}\\
&=&-\frac{1}{2}\sum_{j} e_{j} \cdot T\cdot\nabla^{s}_{e_{j}}\varphi-\frac{1}{2}T\cdot D^{s}(\varphi)\label{slash},
\end{eqnarray}
where $D^{s}$ is the (generalized) Dirac operator induced by $\nabla^{s}$.
\end{lemma}
  \vskip 0.1cm
\noindent {\bf Proof.}  The first formula is an immediate consequence of (\ref{gens1}). For the second description, we use (\ref{clif1}), 
 the definition of $\slashed{D}^{s}$ and   relation  (\ref{symT}). Then, for some arbitrary spinor fields $\varphi, \psi$ we   conclude that
 \begin{eqnarray*}
\langle \slashed{D}^{s}(\varphi), \psi\rangle&=&-\frac{1}{2}\sum_{j}\langle e_{j}\cdot T\cdot\nabla^{s}_{e_{j}}\varphi, \psi\rangle-\frac{1}{2}\sum_{j}\langle T\cdot e_{j}\cdot \nabla^{s}_{e_{j}}\varphi, \psi\rangle\nonumber\\
 &=&-\frac{1}{2}\sum_{j}\langle  e_{j}\cdot T\cdot\nabla^{s}_{e_{j}}\varphi,   \psi\rangle-\frac{1}{2}\sum_{j}\langle e_{j}\cdot \nabla^{s}_{e_{j}}\varphi, T\cdot \psi\rangle\nonumber\\
 &=&-\frac{1}{2}\sum_{j}\langle  e_{j} \cdot T\cdot\nabla^{s}_{e_{j}}\varphi,  \psi\rangle-\frac{1}{2}\langle T\cdot D^{s}(\varphi), \psi\rangle, 
\end{eqnarray*}
which gives rise to (\ref{slash}).  $\blacksquare$ 
%\begin{equation}\label{slash}
%\slashed{D}^{s}(\varphi)=-\frac{1}{2}\sum_{j} e_{j} \cdot T\cdot\nabla^{s}_{e_{j}}\varphi-\frac{1}{2}T\cdot D^{s}(\varphi).
%\end{equation}
\vskip 0.2cm

%\noindent %Unfortunately, without fixing some special condition for the spinor field $\varphi$,   the  first sum in (\ref{slash}) is hard to be treated, since its behaviour depends on the Dirac operator $D^{s}$ and the 4-form $\sigma_{T}$, i.e. 
\noindent Therefore, when the torsion is $\nabla^{c}$-parallel, it is $\sum_{j}e_{j} \cdot T\cdot\nabla^{s}_{e_{j}}\varphi=D^{s}(T\cdot\varphi)-2(1-4s)\sigma_{T}\cdot\varphi={\rm grad}(\gamma)\cdot\varphi+\gamma D^{s}(\varphi)-2(1-4s)\sigma_{T}\cdot\varphi$, where   $\gamma\in {\rm Spec}(T)$ denotes an eigenevalue of $T$, i.e. we assume (without loss of generality) that $\varphi\in\Sigma^{g}_{\gamma}M$ for some real function $\gamma$, not necessarily constant.  Let us begin our investigation with $\nabla^{c}$-parallel spinors, where $\gamma$ is a real constant. 

\begin{prop}\label{parslash}
Consider   a Riemannian spin manifold $(M^{n}, g, T)$ $(n\geq 3)$ with  $\nabla^{c}T=0$, where  $\nabla^{c}:=\nabla^{g}+\frac{1}{2}T$ is the metric connection with skew-torsion $0\neq T\in\Lambda^{3}T^{*}M$. Assume that  $\varphi_{0}\in\Sigma^{g}_{\gamma}M$  is a non-trivial $\nabla^{c}$-parallel spinor and   $\gamma\in {\rm Spec}(T)$ is an eigenvalue of $T$.  Then, $\varphi_{0}$ is an eigenspinor of the operator $\slashed{D}^{s}$ for any $s\in\bb{R}$,   
\[
\slashed{D}^{s}(\varphi_{0})=-\frac{(4s-1)}{4}\Big[T^{2}+2\|T\|^{2}\Big]\cdot \varphi_{0}=-\frac{(4s-1)}{4}\Big[\gamma^{2}+2\|T\|^{2}\Big]\varphi_{0}.
\]
  \end{prop}
\vskip 0.1cm
\noindent {\bf Proof.} % By the definition of $\nabla^{s}, \nabla^{c}$ and $D^{s}$, $D^{c}$, respectively, we have that 
%%\begin{equation*}%\label{nsc}
%\nabla^{s}_{X}\varphi=\nabla^{c}_{X}\varphi+\frac{4s-1}{4}(X\lrcorner T)\cdot\varphi, \quad D^{s}(\varphi)=D^{c}(\varphi)+\frac{3(4s-1)}{4}T\cdot\varphi,
%\end{equation*}
%for any $X\in\Gamma(TM)$ and $\varphi\in\cal{F}^{g}$.  
 % Consequently,   a $\nabla^{c}$-parallel spinor $\varphi_{0}$ satisfies the following two equations
%\begin{equation}\label{usef3}
%\nabla^{s}_{X}\varphi_{0}=\frac{4s-1}{4}(X\lrcorner T)\cdot\varphi_{0}, \quad D^{s}(\varphi_{0})=\frac{3(4s-1)}{4}T\cdot\varphi_{0}.
%\end{equation}
Based on  (\ref{usef3}) and the definition of $\slashed{D}^{s}$, we  see that  any spinor field $\varphi\in\cal{F}^{g}$ satisfies \[
  \slashed{D}^{s}(\varphi)=\slashed{D}^{c}(\varphi)+\frac{(4s-1)}{4}\cal{T}\cdot\varphi=\slashed{D}^{c}(\varphi)+\frac{(4s-1)}{4}\Big[2\sigma_{T}-3\|T\|^{2}\Big]\cdot\varphi,
  \]
 where $\slashed{D}^{c}:=\sum_{j}(e_{j}\lrcorner T)\cdot\nabla^{c}_{e_{j}}$ is the operator associated to $\nabla^{c}$.  Thus, if $\nabla^{c}\varphi_{0}=0$, then  $\slashed{D}^{c}(\varphi_{0})=0$ and the claim immediately follows in combination  with  $\sigma_{T}\cdot\varphi_{0}=\frac{1}{2}(\|T\|^{2}-T^{2})\cdot\varphi_{0}$ (cf. \cite{AF}).  Of course, the same occurs by applying   (\ref{slash}). Indeed,  we rely again on (\ref{usef3}) and compute that
\begin{eqnarray*}
\slashed{D}^{s}(\varphi_{0})&=&-\frac{1}{2}\sum_{j}   e_{j}\cdot T\cdot\nabla^{s}_{e_{j}}\varphi_{0}-\frac{1}{2} T\cdot D^{s}(\varphi_{0})\\
&=&-\frac{4s-1}{8}\sum_{j} e_{j}\cdot T\cdot (e_{j}\lrcorner T)\cdot\varphi_{0} -\frac{3(4s-1)}{8} T^{2}\cdot\varphi_{0}.
 \end{eqnarray*}
However, it is $e_{j}\cdot T=-T\cdot e_{j}-2(e_{j}\lrcorner T)$, hence one can write
\begin{eqnarray*}
\slashed{D}^{s}(\varphi_{0})&=&\frac{4s-1}{8}\sum_{j} T\cdot e_{j}\cdot (e_{j}\lrcorner T)\cdot\varphi_{0}+ \frac{2(4s-1)}{8}\sum_{j}  (e_{j}\lrcorner T)\cdot (e_{j}\lrcorner T)\cdot\varphi_{0} -\frac{3(4s-1)}{8} T^{2}\cdot\varphi_{0}\\
&=&\frac{3(4s-1)}{8} T^{2}\cdot\varphi_{0} +\frac{(4s-1)}{4} (2\sigma_{T}-3\|T\|^{2})\cdot\varphi_{0} -\frac{3(4s-1)}{8} T^{2}\cdot\varphi_{0}\\
&=&\frac{(4s-1)}{4} (2\sigma_{T}-3\|T\|^{2})\cdot\varphi_{0}=\frac{(4s-1)}{4}\cal{T}\cdot\varphi_{0}.
\end{eqnarray*}
Thus the assertion follows by using the relations $\cal{T}=-(T^{2}+2\|T\|^{2})$ and $T^{2}\cdot\varphi_{0}=\gamma^{2}\cdot\varphi_{0}$. $\blacksquare$ %In particular, since $T^{2}=\|T\|^{2}-2\sigma_{T}$ the latter sum is also equal to  $\cal{S}=-T^{2}-2\|T\|^{2}=-(T^{2}+2\|T\|^{2})$. Hence,  
%\begin{equation}\label{fine}
%  \slashed{D}^{s}(\varphi_{0})  =- \frac{(4s-1)}{4}  (T^{2}+2\|T\|^{2})\cdot\varphi_{0} 
%  \end{equation}
%which proves our main claim. Now, if $\nabla^{c}T=0$, then recall from before that the spinor bundle decomposes into the $T$-eigensubundles.  Hence, if our $\nabla^{c}$-parallel spinor $\varphi_{0}$ satisfies in addition $\varphi_{0}\in\Sigma_{\gamma}^{g}M$, then the relation $T\cdot \varphi_{0}=\gamma\varphi_{0}$ implies that  $T^{2}\cdot\varphi_{0}=\gamma^{2}\varphi_{0}$ and by using (\ref{fine}) we conclude. $\blacksquare$ 

\vskip 0.2cm
\noindent The action of the operator $\slashed{D}^{s}$  on Killing spinors and twistor spinors (with torsion or not), with respect to the family $\nabla^{s}$, is known by \cite{Chrysk2}. In particular,  for a non-trivial element $\varphi_{0}\in\ker(\cal{P}^{s})$ for some   $s\in\bb{R}$ and independently of the assumption $\nabla^{c}T=0$, it is not hard to show  that
% by using for example  (\ref{slash}), that
 %\begin{eqnarray*}
%\slashed{D}^{s}(\varphi_{0})&=&-\frac{1}{2}\Big[\sum_{j} e_{j} \cdot T\cdot\nabla^{s}_{e_{j}}\varphi_{0}+T\cdot D^{s}(\varphi_{0})\Big]=\frac{1}{2n}\sum_{j}  e_{j}\cdot T\cdot e_{j}\cdot D^{s}(\varphi_{0})-\frac{1}{2} T\cdot D^{s}(\varphi_{0})\\
%&=&\frac{(n-6)}{2n}T\cdot D^{s}(\varphi_{0})-\frac{1}{2} T\cdot D^{s}(\varphi_{0})=-\frac{3}{n}  T\cdot D^{s}(\varphi_{0}),
%\end{eqnarray*}
%  since $\sum_{j}e_{j}\cdot T\cdot e_{j}=(n-6)T$, see \cite[p.~328]{ABK} for hints. Hence

\begin{prop}\label{twspin} \textnormal{(\cite{Chrysk2})}  Consider   a Riemannian spin manifold $(M^{n}, g, T)$ $(n\geq 3)$  endowed with a non-trivial 3-form $T\in\Lambda^{3}T^{*}M$ and the one-parameter family of metric connections $\nabla^{s}=\nabla^{g}+2sT$. Then, 
any twistor spinor   $\varphi_{0}\in\ker(\cal{P}^{s})$  (with torsion or not), with respect to $\nabla^{s}$ for some $s\in\bb{R}$, satisfies  
\begin{equation}\label{tspin}
\slashed{D}^{s}(\varphi_{0})=-\frac{3}{n}T\cdot D^{s}(\varphi_{0}).
\end{equation}
Moreover,   if $\varphi_{0}\in\cal{K}^{s}(M, g)_{\zeta}$ for  some $s\in\bb{R}\backslash\{0, 1/4\}$ and  $\zeta\neq0$, then 
 $\slashed{D}^{s}(\varphi_{0})=3\zeta T\cdot\varphi_{0}$ and similarly, if  $\varphi_{0}\in\cal{K}^{g}(M, g)_{\kappa}$ for some $\kappa\neq 0$, then $\slashed{D}^{g}(\varphi_{0})=3\kappa T\cdot\varphi_{0}$. 
 \end{prop}
 \begin{corol} \label{corolkil}
 Whenever $\nabla^{c}T=0$, a non-trivial KsT (resp. real Killing spinor) $\varphi_{0}$  induces a non-trivial eigenspinor of $\slashed{D}^{s}$ for the same $s$ (resp. for $s=0$) with eigenvalue $\be:=3\gamma\zeta$, (resp. $\be=3\gamma\kappa$), where    $\gamma\in{\rm Spec}(T)$ is the corresponding $T$-eigenvalue.
\end{corol}
%\vskip 0.1cm
%\noindent {\it Proof.} The relation (\ref{tspin})  follows immediately by combining the definition of a non-trivial  TsT $\varphi_{0}\in\ker(\cal{P}^{s})$, with that of the operator $\slashed{D}^{s}$, see \cite{Chrysk2}.   Of course, any of the equivalent expressions given in Lemma \ref{equivex} yields the same result. $\blacksquare$ 

 %\noindent Examples of special structures satisfying Proposition  \ref{twspin} (with $\be\neq 0$) are  for instance the 6-dimensional nearly K\"ahler manifolds and the 7-dimensional nearly parallel $\G_2$-manifolds. Let us shortly illustrate the case of nearly K\"ahler manifolds. % and in a forthcoming work we shall describe a different class of $\slashed{D}^{s}$-eigenspinors.
\begin{example}\label{nKm}
\textnormal{Consider a 6-dimensional (strict) nearly K\"ahler manifold $(M^{6}, g, J)$, i.e. an almost Hermitan manifold endowed with a   non-integrable  almost complex structure $J$ such that  $(\nabla^{g}_{X}J)X=0$.     By \cite[Thm.~10.1]{FrIv}  it is known that $M^{6}$ admits a (unique)   characteristic connection $\nabla^{c}$ with {\it parallel} skew-torsion, given by $T(X, Y):=(\nabla^{g}_{X}J)JY$. Moreover, there exist two $\nabla^{c}$-parallel spinors $\varphi^{\pm}$ such that $\cal{F}^{g}(\pm 2\|T\|)$, i.e. $T\cdot\varphi^{\pm}=\pm 2\|T\|\cdot\varphi^{\pm}$. Thus, by Proposition \ref{parslash} we get
\begin{equation}\label{lastone}
 \slashed{D}^{s}(\varphi^{\pm})=-\frac{3(4s-1)\|T\|^{2}}{2}\varphi^{\pm}.
 \end{equation}
On the other hand, $\varphi^{\pm}$ are real Killing spinors with $\kappa:=\mp\|T\|/4$,  TsT with torsion for any $s\neq 1/4$, i.e.    $\varphi^{\pm}\in\ke(\cal{P}^{s}|_{\Sigma^{g}_{\pm 2\|T\|}M})$  and  KsT for any $s\neq 0, 1/4$, with Killing number $\zeta:=\mp\frac{(4s-1)}{4}\|T\|$, see \cite[Thm.~4.1]{Chrysk2}.  Therefore, (\ref{lastone}) is deduced also by applying Corollary \ref{corolkil}. }
\end{example}

\vskip 0.1cm
\noindent  Under the condition $\nabla^{c}T=0$,  a kind of converse of Proposition \ref{twspin} reads as follows:
 \begin{prop}\label{endsla}
 Consider a triple $(M^{n}, g, T)$ $(n\geq 3)$ with  $\nabla^{c}T=0$, where  $\nabla^{c}:=\nabla^{g}+\frac{1}{2}T$ is the metric connection with skew-torsion $0\neq T\in\Lambda^{3}T^{*}M$.  Assume that $\varphi_{0}\in\Gamma(\Sigma^{g}_{\gamma})\cap\ke(\cal{P}^{s}):=\ke(\cal{P}^{s}|_{\Sigma_{\gamma}^{g}M})$ is a non-trivial restricted twistor spinor (with torsion or not), for some $s_{0}\in\bb{R}$ and  some non-zero constant   eigenvalue $0\neq\gamma\in{\rm Spec}(T)$. If $\varphi_{0}$ is  a $\slashed{D}^{s_{0}}$-eigenespinor, i.e. $\slashed{D}^{s_{0}}(\varphi_{0})=\be\varphi_{0}$ for some constant eigenvalue $\be$, then
 \begin{equation}\label{interest}
 D^{s_{0}}(\varphi_{0})=\frac{(n-6)\be}{3\gamma}\varphi_{0}+\frac{2(1-4s_{0})}{\gamma}\sigma_{T}\cdot\varphi_{0}.
 \end{equation}
 If $n=6$ or $\be=0$, i.e. $\varphi_{0}\in\ker(\slashed{D}^{s_{0}})$, then $D^{s_{0}}(\varphi_{0})=\frac{2(1-4s_{0})}{\gamma}\sigma_{T}\cdot\varphi_{0}$.
 \end{prop}
 \vskip 0.1cm
 \noindent {\bf Proof.}
 Since $\slashed{D}^{s_{0}}(\varphi_{0})=\be\varphi_{0}$ and $\varphi_{0}\in\ke(\cal{P}^{s_{0}}\big|_{\Sigma_{\gamma}^{g}M})$, the type (\ref{tspin}) reduces to $T\cdot D^{s_{0}}(\varphi_{0})=-\frac{n\be}{3}\varphi_{0}$ and since $\gamma\neq 0$ is a real constant such that $T\cdot\varphi_{0}=\gamma\cdot\varphi_{0}$,  our claim follows by relation (\ref{clasdif1}).  $\blacksquare$
 \begin{corol}
If $\varphi_{0}\in\ke(\cal{P}^{c}\big|_{\Sigma_{\gamma}^{g}M})$ is a non-trivial restricted twistor spinor  with torsion   with respect to $\nabla^{c}$ and $\varphi_{0}$ is $\slashed{D}^{c}$-harmonic, i.e. $\be=0$ and hence $\slashed{D}^{c}(\varphi_{0})=0$, then $D^{c}(\varphi_{0})=0$, in particular $\varphi_{0}$ is $\nabla^{c}$-parallel, i.e.  $\varphi_{0}\in\ker(\nabla^{c})$.
 \end{corol}
 \vskip 0.1cm
 \noindent {\bf Proof.} This follows by Proposition \ref{endsla} in combination with \cite[Lem.~2.2]{Chrysk2}, see also \cite[p.~119]{Chrysk2} for details. $\blacksquare$.
  \vskip 0.2cm
  \noindent We deduce  that the relation $\varphi_{0}\in \ker(\slashed{D}^{c})\cap \ke(\cal{P}^{c}|_{\Sigma_{\gamma}^{g}M})$ for some constant $\gamma\neq 0$, is a very strong condition  which in fact implies the $\nabla^{c}$-parallelism of $\varphi_{0}$, similarly with the condition $\varphi_{0}\in \ker({D}^{c})\cap \ke(\cal{P}^{c}|_{\Sigma_{\gamma}^{g}M})$.  Hence,  in general we avoid to consider this kind of TsT, as in \cite{Chrysk2}.

 \vskip 0.2cm
\noindent  Recall finally by Proposition \ref{parslash} that for a $\nabla^{c}$-parallel spinor field $\varphi_{0}$ the relation $\slashed{D}^{s}(\varphi_{0})=\be\varphi_{0}$ is always verified with $\be=\frac{(4s-1)}{4}\Big[2\sigma_{T}-3\|T\|^{2}\Big]$. Adding  now the extra condition   $\varphi_{0}\in \ke(\cal{P}^{s}|_{\Sigma_{\gamma}^{g}M})$ for some constant $\gamma\neq 0$ and $s\neq 1/4$, then for $3\leq n\leq 8$ we see that relation (\ref{interest}) gives rise to an alternative way to verify that $\varphi_{0}$  is actually a KsT with $\zeta=\frac{3(1-4s)}{4n}$ (see \cite[Thm.~3.7]{Chrysk2}), i.e. $D^{s}(\varphi_{0})=\frac{3(4s-1)\gamma}{4}\cdot\varphi_{0}$  as it should be  according to (\ref{usef3}).
  For such a proof one  may  use the formulas  $\gamma^{2}=\frac{2n}{9-n}\|T\|^{2}$ and     $\sigma_{T}\cdot\varphi_{0}=-\frac{3\gamma^{2}(n-3)}{4n}\varphi_{0}$, given in  \cite[Prop.~3.2]{Chrysk2}.  For $n=6$, relation  (\ref{interest}) is simplified  and  we do not need the explicit form of $\be$.  This case of course applies on nearly K\"ahler manifolds.  For nearly parallel $\G_2$-manifolds $(M^{7}, g, \omega)$ and for the unique $\nabla^{c}$-parallel spinor $\varphi_{0}\in \ke(\cal{P}^{s}|_{\Sigma_{-\sqrt{7}\|T\|}^{g}M})$ one computes  $\be=-\frac{9(4s-1)}{4}\|T\|^{2}$ via Proposition \ref{parslash},  hence   (\ref{interest}) in combination with $\sigma_{T}\cdot\varphi_{0}=-3\|T\|^{2}$   yield the result:
  \[
   D^{s}(\varphi_{0})=-\frac{21(4s-1)}{4\sqrt{7}}\|T\|\cdot\varphi_{0}=\frac{3(4s-1)\gamma}{4}\cdot\varphi_{0},
 \]
 thus $\varphi_{0}\in\cal{K}^{s}(M^{7}, g)_{\zeta}$ with $\zeta=-\frac{3(4s-1)}{4\sqrt{7}}\|T\|$ (cf. \cite[Thm.~4.2]{Chrysk2}). %\noindent In a fortcoming work we will discuss a  different class of $\slashed{D}^{s}$-eigenspinors.}% Observe however, that in other dimensions the equation $\slashed{D}^{s}(\varphi_{0})=\be\varphi_{0}$ fails to induce the condition $\varphi_{0}\in\cal{K}^{s}(M, g)_{\zeta}$ for some  $\varphi_{0}\in \ker({\nabla}^{c})\cap \ke(\cal{P}^{c}\big|_{\Sigma_{\gamma}^{g}M})$, with $s\neq 1/4$. }

 \end{document}